\documentclass[11pt]{amsart}
\usepackage{amscd,amssymb,pstricks,pst-all,fancyvrb,url}
\usepackage[all]{xy}
\newtheorem{theorem}{Theorem}[subsection]
\newtheorem{lemma}[theorem]{Lemma}
\newtheorem{proposition}[theorem]{Proposition}
\newtheorem{corollary}[theorem]{Corollary}

\theoremstyle{definition}
\newtheorem{definition}[theorem]{Definition}
\newtheorem{example}[theorem]{Example}
\newtheorem{remark}[theorem]{Remark}

\newtheorem*{claim}{Claim}
\newtheorem*{ack}{Acknowledgment}
\newtheorem{question}[theorem]{Question}

\newcommand{\A}{{\mathcal A}}

\newcommand{\Mo}{{\mathcal Z}}
\newcommand{\I}{{\mathcal I}}
\renewcommand{\DJ}{{\mathcal{DJ}}}

\newcommand{\Z}{{\mathbb Z}}
\newcommand{\Q}{{\mathbb Q}}
\newcommand{\R}{{\mathbb R}}
\newcommand{\C}{{\mathbb C}}

\newcommand{\CP}{{\mathbb{CP}}}

\renewcommand{\k}{\Bbbk}
\DeclareMathOperator{\Tor}{Tor}
\DeclareMathOperator{\Ext}{Ext}

\DeclareMathOperator{\rank}{rank}
\DeclareMathOperator{\im}{im}

\DeclareMathOperator{\id}{id}

\DeclareMathOperator{\gr}{gr}
\DeclareMathOperator{\Bier}{Bier}
\DeclareMathOperator{\CW}{cw}

\def\g{{\mathfrak g}}

\newrgbcolor{slate}{0.71 0.82 0.84}

\newcommand{\surj}{\twoheadrightarrow}
\newcommand{\inj}{\hookrightarrow}

\def\set#1{{\left\{#1\right\}}}
\def\abs#1{{\left|#1\right|}}

\def\angl#1{{\langle #1\rangle}}
\def\no{{\circ}}
\newcommand{\ones}{{K^{(1)}}}
\DeclareMathOperator{\cs}{{\hbox{\Large \#}}}

\newenvironment{romenum}
{

\begin{enumerate}}{\end{enumerate}}

\begin{document}

\title[Moment-angle complexes, monomial ideals, and 
Massey products]%
{Moment-angle complexes, monomial ideals, and 
Massey products}

\author[G. Denham]{Graham Denham$^1$}
\address{Department of Mathematics, University of Western Ontario,
London, ON  N6A 5B7}
\email{{gdenham@uwo.ca}}
\urladdr{{http://www.math.uwo.ca/\~{}gdenham}}
\thanks{{$^1$}Partially supported by a grant from NSERC of Canada}

\author[A.~I. Suciu]{Alexander~I.~Suciu$^2$}
\address{Department of Mathematics,
Northeastern University,
Boston, MA 02115}
\email{{a.suciu@neu.edu}}
\urladdr{{http://www.math.neu.edu/\~{}suciu}}
\thanks{$^2$Partially supported by NSF grant DMS-0311142}



\dedicatory{To Robert MacPherson on the occasion of his sixtieth birthday}

\begin{abstract}
Associated to every finite simplicial complex $K$ 
there is a ``moment-angle" finite CW-complex, 
$\Mo_K$;   if $K$ is a triangulation of a sphere, 
$\Mo_K$ is a smooth, compact manifold.  Building 
on work of Buchstaber, Panov, and Baskakov, we 
study the cohomology ring, the homotopy groups, 
and the triple Massey products of a moment-angle 
complex, relating these topological invariants 
to the algebraic combinatorics of the underlying 
simplicial complex.   Applications to the study 
of non-formal manifolds and subspace 
arrangements are given.

\end{abstract}

\maketitle

\setcounter{tocdepth}{1}
\tableofcontents

\section{Introduction}
\label{sec:intro}

\subsection{Moment-angle complexes}
\label{intro:mom}

A construction due to Davis and Januszkiewicz \cite{DJ91} 
and studied in detail by Buchstaber and Panov \cite{BP00} 
associates to every simplicial complex $K$ on $n$ vertices 
a finite cellular complex $\Mo_K$, endowed with a natural 
action by the $n$-torus, and whose orbit space is the 
cone over $K$.   

The unit disk $D^2\subset \C$ has a 
natural cell structure, in which the boundary circle $S^1$ is the 
$1$-skeleton.  The $n$-polydisk $(D^2)^{\times n}$ 
inherits a product cell decomposition, with the boundary 
$n$-torus, $T^{n}=(S^1)^{\times n}$, as a subcomplex.  
For each subset $\sigma\subseteq [n]:=\set{1,\ldots,n}$, let 
$B_{\sigma}=\{ z \in (D^2)^{\times n} \mid \abs{z_i}=1\ 
\text{for $i\notin\sigma$} \}$. Taking the union of the 
subcomplexes $B_{\sigma}$, indexed by the simplices 
of $K$, yields the {\em moment-angle complex} $\Mo_K$.  
This cellular complex is always $2$-connected and has
dimension $n+\ell+1$, where $\ell$ is the dimension of $K$.

It turns out that the algebraic topology of a moment-angle 
complex $\Mo_K$, as embodied, for example, in the 
cohomology ring, the homotopy groups, and the Massey 
products, is intimately related to the combinatorics 
of the underlying simplicial complex $K$. 

\subsection{Cohomology ring and homotopy groups}
\label{intro:homotopy}

Fix a coefficient field $\k$ of characteristic $0$.  
Let $S=\k[x_1,\dots,x_n]$ be the polynomial ring 
with variables in degree $2$, and let $I$ 
be the ideal generated by all monomials 
corresponding to non-faces of $K$.  
As shown by Buchstaber and Panov \cite{BP00}, the 
cohomology ring of the moment-angle complex, 
$H^*(\Mo_K,\k)$, is isomorphic to $\Tor^S(S/I,\k)$, 
the Tor-algebra of the Stanley-Reisner ring $S/I$.  

Here, we compute the ranks of the homotopy 
groups of $\Mo_K$ in terms of the homological 
algebra of $S/I$. A crucial ingredient is provided 
by the fibration $\Mo_K \to \DJ(K) \to BT^n$, where  
$\DJ(K)$ is the Davis-Januszkiewicz space associated 
to $K$.  The formality of this space---a result of 
Franz~\cite{Fr03} for smooth toric fans, and Notbohm 
and Ray~\cite{NR05} in general---leads to the collapse of the 
Eilenberg-Moore spectral sequence for the path-fibration 
of $\DJ(K)$, thereby allowing us to compute 
$\pi_*(\Mo_K) \otimes \k$ in terms of $\Tor^{S/I}(\k,\k)$, 
the dual of the Yoneda algebra of $S/I$.  

The answer turns out to be particularly nice when $K$ 
is a flag complex, in which case $S/I$ is a Koszul algebra.  
Writing $h(H,s,t)$ for the bigraded Hilbert series of 
$H=\Tor^S(S/I,\k)$, and making use of Koszul duality, 
we prove (Theorem \ref{thm:flag}):
\begin{equation*}
\label{eq:lag}
\prod_{r=1}^{\infty} 
(1-(-t)^{r})^{(-1)^r \rank \pi_{r+1}(\Mo_K)}=
h(H,i\sqrt{t},-i\sqrt{t}).
\end{equation*}
This formula provides an effective method for computing 
the rational homotopy groups of such moment-angle complexes, 
at least in low degrees; see Examples \ref{ex:polygons} 
and \ref{ex:indec}.

\subsection{Cellular cochains and Massey products}
\label{intro:cellsmassey}

A key observation of Buchstaber and Panov \cite{BP00} is that 
the cellular cochain complex $C^*_{\CW}(\Mo_K,\k)$ comes 
endowed with a multiplication which makes it into a 
commutative differential (bi)graded algebra, quasi-isomorphic 
to the Koszul complex of $S/I$.  An explicit formula for 
the cochain multiplication, in terms of pairings between 
full subcomplexes, can be found in \cite{Ba02}.  

Using this approach, Baskakov \cite{Ba03} constructs 
simplicial complexes $K$ for which the cellular cochain 
algebra $C^*_{\CW}(\Mo_K,\k)$ has non-vanishing Massey 
triple products.  In turn, this implies the non-formality of 
$\Mo_K$; see \S\ref{sec:cellmassey} for more details. 

We embark here on a more systematic study of Massey 
products in moment-angle complexes.   As a start, we 
characterize those simplicial complexes $K$ for which 
$\Mo_K$ has a non-trivial Massey product in lowest 
possible degree, purely in combinatorial terms.  
In Theorem \ref{th:char}, we prove: 
there exist $\alpha, \beta, \gamma \in H^3(\Mo_K,\k)$ 
such that  $\angl{\alpha,\beta,\gamma}\ne 0$ 
precisely when the $1$-skeleton of $K$ has 
an induced subgraph isomorphic to one 
of the five ``obstruction" graphs listed  in 
Figure~\ref{fig:exclude}. (The smallest such 
simplicial complex is depicted in Figure~\ref{fig:sixvertex}.)
Moreover, all Massey products arising in this fashion 
are decomposable.  

\begin{figure}
\begin{pspicture}(-5.0,-1.2)(5,3.3)
\definecolor{wheat}{rgb}{.96,.87,.7}
\definecolor{lightblue}{rgb}{.68,.85,.9}
\cnode*[linecolor=black](-2,-1){0.06}{1}
\cnode*[linecolor=black](2,2){0.06}{2}
\cnode*[linecolor=black](0,0){0.06}{3}
\cnode*[linecolor=black](-2,3){0.06}{4}
\cnode*[linecolor=black](2,0){0.06}{5}
\cnode*[linecolor=black](0,2){0.06}{6}
\rput(-2.2,-1){1}
\rput(2.25,2){2}
\rput(0.2,0.2){3}
\rput(-2.2,3){4}
\rput(2.25,0){5}
\rput(0.2,1.8){6}
\pspolygon[fillstyle=solid, fillcolor=slate](0,2)(-2,3)(2,2)
\pspolygon[fillstyle=solid, fillcolor=wheat](0,2)(-2,3)(-2,-1)
\pspolygon[fillstyle=solid, fillcolor=slate](0,2)(-2,-1)(0,0)
\pspolygon[fillstyle=solid, fillcolor=wheat](-2,-1)(0,0)(2,0)
\ncline{2}{5}
\end{pspicture}
\caption{\textsf{A simplicial complex $K$ for 
which $\Mo_K$ is non-formal}}
\label{fig:sixvertex}
\end{figure}

On the other hand, we exhibit in \S\ref{subsec=indecomp} 
a simplicial complex $K$ for which $\Mo_K$ 
carries an indecomposable triple Massey product;  
such a product manifests itself as a non-zero 
differential in the $E_2$ term of the Eilenberg-Moore 
spectral sequence for $\Mo_K$.

\subsection{Moment-angle manifolds}
\label{intro:triangulations}

Now suppose $K$ is an $n$-vertex triangulation of the sphere 
$S^{\ell}$.  Then, as shown by Buchstaber and Panov \cite{BP00}, 
the moment-angle complex $\Mo_K$ is a smooth, compact, 
$2$-connected manifold of dimension $n+\ell+1$.  
Some of these manifolds can be described in simple terms. 
For example, if the triangulation $K$ is dual to a polytope 
obtained from the hypercube by cutting corners,  then $\Mo_K$ 
is a connected sum of products of spheres.   

In general, though, moment-angle manifolds can exhibit quite 
a complicated structure, both from the point of view of 
their cohomology ring and of their Massey products. 
We illustrate this point with an infinite family 
of triangulations of $S^2$ for which the corresponding 
manifolds have non-trivial triple Massey products.  
The simplest such manifold is constructed as follows: 
start with the square $\square$; form the deleted join of 
$\square$ with its Alexander dual, $\times$, to obtain 
$K=\Bier(\square)$, an $8$-vertex Bier triangulation 
of $S^2$; the resulting moment-angle complex, $\Mo_K$, 
is a non-formal, $11$-dimensional manifold.

Asymptotically, almost all triangulations of $S^2$ yield 
non-formal moment-angle manifolds: see Theorem \ref{th:asym} 
for a precise statement.   We detect this non-formality by 
means of decomposable Massey products. Using the Bier sphere 
construction, we exhibit in \S\ref{subsec:bier} a $16$-vertex 
triangulation $K$ of $S^6$ for which $\Mo_K$ has an
{\em indecomposable}\/  Massey product. 

\subsection{Compact, complex, non-K\"ahler manifolds}
\label{intro:lmv}

Classical constructions of complex manifolds, due 
to Hopf and Calabi--Eckmann were generalized 
in recent years by L\'opez de Medrano--Verjovsky 
\cite{LV97} and Meersseman \cite{Me00}.  These authors 
define a large class of compact complex manifolds 
admitting no K\"ahler structure. A complex manifold 
$N$ arises via the LVM construction as the leaf space 
of a foliation of $\CP^{n-1}$, given by a suitable 
linear action of $\C^m$ on $\C^n$, with $n>2m$. 
If $n=2m+1$, then $N$ is a complex torus; otherwise, 
$N$ is non-symplectic, and thus non-K\"ahler.  

It turns out that the LVM construction is very much 
related to that of the moment-angle complexes.  
We find some parallels and applications both ways, 
especially in connection with recent work of 
Bosio and Meersseman \cite{BM04}. In particular, 
we explain the absence of K\"ahler structure for certain 
LVM manifolds via the presence of non-vanishing Massey 
products in the corresponding moment-angle complex.

\subsection{Subspace arrangements}
\label{intro:subs}

Much is known about the relationship between the combinatorics 
of a complex subspace arrangement (as encoded in its ranked 
intersection lattice) and the topology of its complement. 
In \cite{GM88}, Goresky and MacPherson gave a formula 
for the Betti numbers of the complement; the cup-products 
in cohomology were computed by Deligne-Goresky-MacPherson 
\cite{DGM00} and de Longueville-Schultz \cite{dLS01}.  

Determining the homotopy type of the complement remains 
a challenging problem.  If the intersection lattice is geometric, 
then, as shown by Feichtner and Yuzvinsky \cite{FY05}, 
the complement is formal; thus, its rational homotopy type 
is determined by the (ranked) intersection lattice.  
Using our approach, we see that the complement of a subspace 
arrangement  is \emph{not}\/ formal in general:  non-vanishing 
triple Massey products are not detected by the cohomology ring.  
The simplest such example is the arrangement of coordinate 
subspaces in $\C^6$ determined by the above simplicial 
complex: $\A =\{H_{1} ,\dots , H_{5}\}$, 
where $H_{i}=\{z \in \C^6 \mid z_i=z_{i+1}=0\}$.

\section{The moment-angle functor}
\label{sec:mom}

\subsection{Generalized moment-angle complexes}
\label{subsec:mogen}

We start with a generalization of the notion of 
moment-angle complex, due to Strickland \cite{Str99}, 
cf.~\cite{BP00, Ba02, Pa05}.  

\begin{definition}
\label{def:zkx}
Let $X$ be a space, and $A\subset X$ a non-empty 
subspace.  Given a simplicial complex $K$ on vertex set 
$[n]$, define $\Mo_K(X,A)$ to be 
the following subspace of the cartesian product $X^{\times n}$:
\begin{equation}
\label{eq:zkx}
\Mo_K(X,A)=\bigcup_{\sigma\in K}  (X,A)^{\sigma}, 
\end{equation}
where $(X,A)^{\sigma}= \{ x \in X^{\times n} \mid x_i \in A \text{ if } 
i\notin \sigma\}$.
\end{definition}

In other words, $\Mo_K(X,A)$ is the colimit of the diagram 
of spaces $\{(X,A)^{\sigma}\}$, indexed by the category 
whose objects are the simplices of $K$ (including the 
empty simplex), and whose morphisms are the inclusions 
between those simplices; see \cite{PRV04, NR05}.

Note that if $A=X$, or $K=\Delta^{n-1}$ is a simplex, then 
$\Mo_K(X,A)=X^{\times n}$; at the other extreme, 
$\Mo_{\emptyset}(X,A)=A^{\times n}$.  From the 
definition, if $L$ is a subcomplex of $K$ on the 
same vertex set, then $\Mo_L(X,A)$ is a subspace 
of $\Mo_K(X,A)$; in particular, 
$A^{\times n}\subset\Mo_K(X,A)\subset X^{\times n}$.  

If $X$ is a finite-type CW-complex, and $A$ 
is a subcomplex, then $\Mo_K(X,A)$ is 
a subcomplex of the product complex $X^{\times n}$; 
if $X$ is a finite CW-complex, $\Mo_K(X,A)$ 
is also a finite CW-complex.  

\begin{example}
\label{ex:mompt}
If $(X,*)$ is a pointed space, we will simply write 
$\Mo_K(X):=\Mo_K(X,*)$.  Particularly noteworthy  
are the following examples.

\begin{enumerate}
\item
If $K$ is a discrete set of $n$ points, then 
$\Mo_K(X)=\bigvee^n X$, the wedge of $n$ copies 
of $X$. 

\item If $K=\partial \Delta^{n-1}$ is the boundary of 
a simplex, then $\Mo_K(X)$ is the fat wedge of $n$ 
copies of $X$.

\item 
If $\Delta(\Gamma)$ is the flag complex of a simple 
graph $\Gamma$, then $\Mo_{\Delta(\Gamma)}(S^1)=
K(G_{\Gamma},1)$, an Eilenberg-MacLane space 
for the right-angled Artin group $G_{\Gamma}$; 
see \cite{CD95, MV95}.

\item  
If $X$ is the classifying space $BS^1=\CP^{\infty}$, 
then $\Mo_K(BS^1)=\DJ(K)$, 
the Davis-Januszkiewicz space associated to $K$; 
see \cite{DJ91, BP00},  and also \S\ref{subs:fibrations}.

\end{enumerate} 
\end{example}

\begin{example}
\label{ex:mom}
In the case when $(X,A)=(D^2,S^1)$, we obtain  
the usual moment-angle complex, $\Mo_K=\Mo_K(D^2,S^1)$, 
which is our main object of study.  Note that 
$\dim \Mo_K = \dim K + n + 1$.  Here are some 
samples, to get started.

\begin{enumerate}
\item If $K=\partial \Delta^{n-1}$, 
then  $\Mo_K=S^{2n-1}$.  

\item If $K$ a discrete set of $n$ points, then 
$\Mo_{K} \simeq \bigvee_{k=2}^{n} (k-1)\binom{n}{k} S^{k+1}$; 
see \cite{GT04}. 
\end{enumerate}
\end{example}

The moment-angle construction 
behaves nicely with respect to joins of simplicial complexes.  

\begin{lemma}
\label{lem:joins}
For any simplicial complexes $K$, $K'$ and pair of 
spaces $(X,A)$, there is a homeomorphism 
\begin{equation*}
\Mo_{K*K'}(X,A) \xrightarrow{\,\cong\,} \Mo_{K}(X,A)\times\Mo_{K'}(X,A).
\end{equation*}
 \end{lemma}

\begin{proof}
View the respective moment-angle complexes as 
subspaces $\Mo_{K}(X,A)\subset X^{\times n}$, 
$\Mo_{K'}(X,A)\subset X^{\times n'}$, and 
$\Mo_{K*K'}(X,A)\subset X^{\times (n+n')}$.  
The natural identification map $f\colon X^{\times(n+n')}  
\xrightarrow{\,\cong\,}X^{\times n}\times X^{\times n'}$ 
restricts to the desired homeomorphism.
\end{proof}

\begin{remark}
\label{rem:extend}
One may slightly generalize Definition \ref{def:zkx}, 
and allow the vertex set $V$ of $K$ to be a (possibly 
proper) subset of $[n]$.  The resulting subspace of 
$X^{\times n}$, call it $\Mo_{K,n}(X,A)$, is then 
homeomorphic to $\Mo_{K}(X,A) \times A^{\times (n-\abs{V})}$. 
For example, $\Mo_{\emptyset,n}(X,A)=A^{\times n}$, and 
$\Mo_{K,n}(D^2,S^1) = \Mo_{K}(D^2,S^1) \times T^{n-\abs{V}}$, where
we recall $T^m$ is the $m$-torus.  We will 
not pursue this generality, except briefly in \S\ref{subsec:complex}. 
\end{remark}

\subsection{Naturality}
\label{subsec:natural}
The $\Mo$ construction enjoys functoriality properties in both 
arguments, which we now summarize.  

\begin{lemma}
\label{lem:zkmap}
Let $f\colon (X,A)\to (Y,B)$ be a map of pairs.  Then 
$f^{\times n} \colon X^{\times n} \to Y^{\times n}$ 
restricts to a map 
\begin{equation}
\label{eq:zkmap}
\Mo_K(f)\colon \Mo_K(X,A) \to \Mo_K(Y,B).
\end{equation}
Moreover, $\Mo_K(g\circ f)=\Mo_K(g) \circ \Mo_K(f)$.  
And, if $f$ is a cellular map between CW pairs, then 
$\Mo_K(f)$ is also a cellular map.  
\end{lemma}

Now suppose $F\colon (X,A)\times I \to (Y,B)$ is a relative homotopy. 
Then the product homotopy  $(X^{\times n},A^{\times n})\times I 
\to (Y^{\times n},B^{\times n})$ restricts to a homotopy 
$\Mo_K(X,A)\times I \to \Mo_K(Y,B)$.  Thus, $\Mo_K$ is 
a homotopy functor: the homotopy type of $\Mo_K(X,A)$ 
depends only on the relative homotopy type of the pair $(X,A)$. 
Also note the following: if $(X,A)$ deform-retracts onto $(Y,B)$, 
then $\Mo_K(X,A)$ deform-retracts onto $\Mo_K(Y,B)$.

The following observation is due to Strickland \cite{Str99}.

\begin{lemma}
\label{lem:zklmono}
Let $X$ be a commutative topological monoid, and 
$A$ a sub-monoid.  Suppose $f\colon K\to L$ is  
simplicial map between simplicial complexes on 
vertex sets $[n]$ and $[m]$. Then, the map 
$\tilde{f}\colon X^{\times n} \to X^{\times m}$ 
defined by $\tilde{f}(x)_j = \prod_{i: f(i)=j} x_i$ 
restricts to a map 
\begin{equation*}
\label{eq:zkl}
\Mo_f(X,A) \colon \Mo_K(X,A) \to \Mo_L(X,A). 
\end{equation*}  
\end{lemma}

Clearly, $\Mo_{g\circ f}(X,A)=\Mo_g(X,A) \circ \Mo_f(X,A)$.  
Moreover, if $(X,A)$ is a CW pair, then $\Mo_f(X,A)$ is a 
cellular map.

The Lemma applies to the pair $(D^2,S^1)$, 
with monoid structure defined by complex multiplication.  
This yields a moment-angle functor from finite simplicial 
complexes  $K$  and simplicial maps $f\colon K\to L$ 
to finite CW-complexes $\Mo_K$ and cellular maps  
$\Mo_f\colon \Mo_K\to \Mo_L$.

There is a second functorial property that holds for all pairs $(X,A)$,
but only certain simplicial maps.  Recall that if $K$ is a subcomplex
of $L$, then $K$ is a {\em full subcomplex} if every simplex of $L$
on vertices of $K$ is also a simplex of $K$.

\begin{lemma}
\label{lem:full}
Let $(X,A)$ be a pair of spaces, and  $*\in A$ a basepoint.  
Suppose $f\colon K\inj L$ 
is the inclusion of a full subcomplex, sending vertex 
set $[n]$ to $[m]$.  Then:
\begin{romenum}
\item 
The canonical projection $X^{\times m} \surj X^{\times n}$ 
restricts to a surjective map
$\Mo^f(X,A) \colon$ $\Mo_L(X,A) \surj \Mo_K(X,A).$
\item  
The map $\hat{f}\colon X^{\times n} \inj X^{\times m}$, 
defined by $\hat{f}(x)_j = x_i$ if $f(i)=j$, and 
$\hat{f}(x)_j = *$, otherwise, restricts to an injective map  
$\Mo_{f,*} (X,A) \colon \Mo_K(X,A) \inj \Mo_L(X,A)$.  
\end{romenum}
Moreover, $\Mo^f(X,A) \circ \Mo_{f,*}(X,A) =\id$.
\end{lemma}

In the case when $X$ is a commutative monoid, 
$A$ is a sub-monoid, and $*\in A$ is the unit, we have 
$\Mo_{f,*}(X,A)=\Mo_f(X,A)$.

\subsection{The homotopy fibre interpretation}
\label{subs:fibrations}

An important property of the functor $\Mo_K$ is that it takes 
certain relative fibrations to fibrations.  We make this observation 
more precise in the following Lemma.  As an application, we 
obtain two useful fibrations (Lemmas~\ref{lem:bg} and \ref{lem:Gfibre}) 
generalizing those of Buchstaber and Panov.

\begin{lemma}
\label{lem:fibrations}
Let $p\colon (E,E') \to (B,B')$ be a map of pairs, such that 
both $p\colon E\to B$ and $p|_{E'}\colon E' \to B'$ are 
fibrations, with fibres $F$ and $F'$, respectively.  
Suppose that either $F=F'$ or $B=B'$.
Then the product fibration, $p^{\times n}  \colon  
E^{\times n}\to B^{\times n}$, 
restricts to a fibration 
\begin{equation}
\label{eq:zkfib}
\xymatrixcolsep{24pt}
\xymatrix{\Mo_K(F,F') \ar[r]& \Mo_K(E,E')\ar^{\Mo_K(p)}[r] 
& \Mo_K(B,B') }.
\end{equation}
Moreover,  if $(F,F')\to (E,E') \to (B,B')$ is a relative bundle 
(with structure group $G$), and either $F=F'$ or $B=B'$, 
then \eqref{eq:zkfib} is also a bundle 
(with structure group $G^{\times n}$).  
\end{lemma}

\begin{proof}
Suppose first that $B=B'$; then $\Mo_K(B,B')=B^{\times n}$.  The
fibre of $p^{\times n}|_{\Mo_K(E,E')}$ equals
\begin{eqnarray*}
F^{\times n}\cap \Mo_K(E,E') 
&=& \bigcup_{\sigma\in K} F^{\times n}\cap (E,E')^\sigma \\
&=&\bigcup_{\sigma\in K} (F,F')^\sigma
\qquad\qquad\hbox{since $F\cap E'=F'$}\\
&=&\Mo_K(F,F').
\end{eqnarray*}
The homotopy lifting property (or the bundle map property) 
for $\Mo_K(p)$ follows from those of $p$ and $p|_{E'}$.  

Similarly, if $F=F'$, it is straightforward to check that
$\Mo_K(E,E')$ is the pullback of $E^{\times n}$ along 
the inclusion of $\Mo_K(B,B')$ into $B^{\times n}$. 
\end{proof}

For a topological group $G$, denote by $G\to EG \to BG$ 
the universal principal $G$-bundle, with $EG$ a contractible 
space endowed with a free $G$-action, and $BG$ the orbit space. 
We may view $G$ as a subspace of $EG$ (the orbit of a 
basepoint $*$).  If $X$ is a $G$-space, denote by 
$EG\times_{G} X$ the Borel construction, and by 
$X\to EG\times_{G} X \to BG$ the associated bundle. 

Now suppose $A$ is a $G$-invariant subspace of $X$.  
The $G$-action on the pair $(X,A)$ extends canonically to an action 
of the product group $G^{\times n}$ on the pair 
$(X^{\times n},A^{\times n})$. It is readily checked 
that this $G^{\times n}$-action preserves the 
subspace $\Mo_K(X,A)\subset X^{\times n}$. 

\begin{lemma}
\label{lem:bg}
With notation as above, we have:
\begin{enumerate}
\item  \label{pt1}
$E G^{\times n} \times_{G^{\times n}} \Mo_K(EG, G)\simeq 
 \Mo_K(BG)$. 
 
\item \label{pt2}
The homotopy fibre of the inclusion 
$\Mo_K(BG) \inj BG^{\times n}$ is $\Mo_K(EG,G)$.  
In other words, we have a fibration sequence  
$\Mo_K(EG,G)\to \Mo_K(BG)\to BG^{\times n}$. 
\end{enumerate}
\end{lemma}

\begin{proof}
The free action of $G$ on $EG$ gives an 
associated bundle $EG\to EG\times_G EG\to BG$.  Restricting 
each fibre $EG$ to the subspace $G$, we obtain a 
sub-bundle $G\to EG\times *\to BG$, isomorphic to the 
classifying bundle.  From Lemma~\ref{lem:fibrations}, 
since the base spaces agree, we obtain the fibre bundle 
sequence in the top row of the following diagram:
\begin{equation}
\label{eq:big diagram}
\xymatrix@R18pt@C3pt{
\Mo_K(EG,G)\ar^(.38){j}[rr] && \Mo_K(EG\times_G EG,EG\times *) 
\ar[rr]&& \Mo_K(BG,BG)\\
&G^{\times n} \ar^(.4){i_0}[rr] \ar_{j_0}[dl] && EG^{\times n} \ar[dl] 
\ar_(.4){i}[ul]& \\
\Mo_K(EG,G)\ar[rr] \ar@{=}[uu] &&
E G^{\times n} \times_{G^{\times n}} \Mo_K(EG, G)
\ar[rr]\ar_(.4){\phi}@{-->}[uu]    && BG^{\times n}  \ar@{=}[uu] \\ 
}
\end{equation}

The bottom row is the associated bundle to 
$G^{\times n} \to EG^{\times n} \to BG^{\times n}$ 
and the natural action of $G^{\times n}$ on 
$\Mo_K(EG,G)\subset E G^{\times n}$. 
The total space of this bundle is the pushout of 
$EG^{\times n}$ and $\Mo_K(EG,G)$ along the 
natural inclusions $i_0$ and $j_0$ of $G^{\times n}$ into each. 

Identifying $EG^{\times n}$ with $(EG\times *)^{\times n}$ 
yields an inclusion $i$ into $\Mo_K(EG\times_G EG,EG\times *)$.  
We also have the inclusion $j$ of $\Mo_K(EG,G)$ into 
that space, as the typical fiber of the top bundle. 
It is readily seen that $i$ and $j$ agree on the common subspace 
$G^{\times n}$; therefore, there is a map $\phi$ as indicated 
in the diagram.  Since both $i$ and $j$ are 
$G^{\times n}$-equivariant, $\phi$ is a bundle map.  
Since the restriction of $\phi$ to each fibre is the identity, 
we conclude that $\phi$ is a homeomorphism.  

Now note that 
$(EG\times_G EG,  EG\times * ) \simeq (BG, *)$, and 
so $\Mo_K(EG\times_G EG,  EG\times * ) \simeq \Mo_K(BG)$, 
by the remark following Lemma~\ref{lem:zkmap}. 
Putting things together finishes the proof of claim \eqref{pt1}. 
Claim \eqref{pt2} follows at once. 
\end{proof}

\begin{lemma}\label{lem:Gfibre}
There is a fibre bundle sequence
\begin{equation*}
G^{\times n}\to \Mo_K(EG,G)\to \Mo_K(BG).
\end{equation*}
\end{lemma}
\begin{proof}
Take the universal bundle $G\to EG\to BG$ and restrict the
base to the basepoint $*$ to obtain a relative $G$-bundle 
$(G,G)\to (EG,G)\to (BG,*)$.  Then apply Lemma~\ref{lem:fibrations}, 
this time where the two fibres agree.
\end{proof}

For the circle $S^1$, the universal bundle is the Hopf fibration 
$S^1\to S^{\infty} \to \CP^{\infty}$.  Thus, a classifying space for 
the torus $T^n=(S^1)^{\times n}$ is $BT^n=(\CP^{\infty})^{\times n}$. 
The circle acts on the pair $(D^2,S^1)$ by rotation.  As noted 
above, this $S^1$-action extends to a $T^n$-action 
on the moment-angle complex $\Mo_K=\Mo_K(D^2,S^1)$.  
Finally, observe that $(ES^1, S^1)\simeq (D^2,S^1)$; 
thus, $\Mo_K(ES^1, S^1)\simeq \Mo_K$. 

The Davis-Januszkiewicz space associated to the simplicial 
complex $K$ is, by definition, the Borel construction 
on the moment-angle  complex $\Mo_K$, viewed as 
a $T^n$-space: 
\begin{equation}
\DJ(K)=E T^n\times_{T^n}\Mo_K.
\end{equation}
We thus have a fibre bundle sequence $\Mo_K\to \DJ(K)\to BT^n$. 
It follows that the homotopy fibre of the inclusion $\Mo_K \inj \DJ(K)$ 
is $T^n$, a particular case of Lemma \ref{lem:Gfibre}.  

From Lemma \ref{lem:bg}, we recover 
the following result of Buchstaber and Panov \cite{BP00}.  

\begin{corollary}
\label{cor:djbp}
The following hold:
\begin{enumerate}
\item  $\DJ(K)\simeq \Mo_K(BS^1)$
\item The homotopy fibre of the inclusion 
$\Mo_K(BS^1)\inj (BS^1)^{\times n}$ is $\Mo_K$. 
\item We have a fibration sequence 
$\Mo_K\to \Mo_K(BS^1) \to BT^n$. 
\end{enumerate}
\end{corollary}

\section{Cohomology and the Stanley-Reisner ring}
\label{section:coho}

In this section, we outline Buchstaber and Panov's 
computation of the cohomology ring of a moment-angle 
complex $\Mo_K$.  We work over a fixed coefficient ring
which is, by default, the integers $\Z$.  We will use $\k$ to 
denote an arbitrary field of characteristic zero.

\subsection{Cohomology ring of $\Mo_K(X)$}
\label{subs:sr}
 
Consider a pointed space $(X,*)$. Each inclusion of simplices
$\sigma\subseteq\tau$ in $K$ gives rise to an inclusion $(X,*)^\sigma
\hookrightarrow (X,*)^\tau$.  Such an inclusion has a left inverse, 
obtained by mapping $X$ to $*$ in coordinates from 
$\tau\setminus \sigma$, and by mapping $X$ to $X$ identically 
elsewhere.  It follows that the induced maps in cohomology, 
$H^*((X,*)^\tau,\Z) \rightarrow H^*((X,*)^\sigma,\Z)$, 
are split surjections.  Due to this splitting, the Mayer-Vietoris 
spectral sequence degenerates at its $E_2$ term, yielding an 
isomorphism 
\begin{equation}
\label{eq:invlim}
H^*(\Mo_K(X),\Z)\cong \varprojlim_{\sigma \in K} 
H^*((X,*)^\sigma,\Z).
\end{equation}

\begin{lemma}
\label{lem:X*quot}
For any pair $(X,*)$, the inclusion 
$j\colon \Mo_K(X,*)\hookrightarrow\Mo_K(X,X)=
X^{\times n}$ induces a surjection of rings
\[
j^*\colon H^*(X^{\times n},\Z)\twoheadrightarrow H^*(\Mo_K(X),\Z).
\]
\end{lemma}
\begin{proof}
We use the same splitting argument.  For each $\sigma\in K$, the
inclusion $(X,*)^\sigma\hookrightarrow (X,X)^\sigma$ has a left
inverse, so the inclusion $j:\Mo_K(X,*)\hookrightarrow\Mo_K(X,X)$
does too.  Therefore, $j$ induces a split surjection in cohomology.
\end{proof}

If $X$ is a CW-complex, the inclusion $j$ above is easily seen to be 
cellular.  Its image consists of cells $e$ for which there exists
a simplex $\sigma\in K$ with $e_i=*$ for all $i\not\in\sigma$. 
If further, $X$ has a minimal cell structure---i.e., if the number 
of $p$-cells of $X$ equals $b_p(X)$, for all $p$---more can 
be said: the kernel of $j^*$ is spanned (additively) by the 
dual basis to the complement of the image of $j$; that is, 
those $e^*$ for which, for some non-face $\sigma$, 
$(e^*)_i\neq1$ for all $i\in\sigma$.

Two important cases are $X=S^1$ and $X=BS^1$; see 
Section~\ref{subsec:mogen}.  In either case, associate 
with a simplex $\sigma=\set{i_1, \ldots,i_k}$ a square-free 
monomial $x_\sigma=x_{i_1}x_{i_2}\cdots x_{i_k}$.

\begin{theorem}[\cite{DJ91, BP00, KR}]
\label{thm:cohoring}
If $K$ is a simplicial complex on $n$ vertices,
\begin{enumerate}
\item $H^*(\Mo_K(S^1),\Z)\cong\bigwedge [x_1,\ldots,x_n]/J_K$, and
\item $H^*(\Mo_K(BS^1),\Z)\cong\Z[x_1,\ldots,x_n]/I_K$,
\end{enumerate}
where $J_K$ and $I_K$ are the ideals in the exterior and 
polynomial algebras, respectively, generated by all monomials 
$x_\sigma$ for which $\sigma$ is not a face of $K$.
\end{theorem}

\begin{proof}
Note that $\C P^\infty=BS^1$ and $S^1$ both have a minimal cell 
structure. We have $H^*((BS^1)^{\times n},\Z)\cong\Z[x_1,\ldots,x_n]$, 
where the generator $x_i$ is dual to the $2$-cell in the $i$th 
coordinate, and $H^*((S^1)^{\times n},\Z)$ is an exterior algebra 
with degree-$1$ generators.

By Lemma~\ref{lem:X*quot}, then, it is enough to notice that the dual
basis to the cells in $\Mo_K(BS^1)$ that are not in $(BS^1)^{\times n}$
consists exactly of the monomials in $I_K$, and similarly for $S^1$.
\end{proof}

Let $E$ and $S$ denote the exterior and polynomial rings above, 
with generators in degree $1$ and $2$, respectively.
That is, the cohomology of the Davis-Januszkiewicz space 
$\DJ(K)= \Mo_K(BS^1)$ is simply the Stanley-Reisner ring of $K$ 
(see \cite{DJ91,BP00,Pa05}):
\begin{equation}
\label{eq:sr}
H^*(\DJ(K),\Z)=S/ I_K.
\end{equation}

When $K=\Delta(\Gamma)$ is the flag complex of 
a graph $\Gamma$, recall from Section \ref{subsec:mogen} 
that $\Mo_K(S^1)$ is a $K(G_\Gamma,1)$ space for 
the right-angled Artin group $G_\Gamma$.  Thus, 
the cohomology of $K(G_\Gamma,1)$ is simply the 
exterior Stanley-Reisner ring of $K$ (see \cite{KR}):
\begin{equation}
\label{eq:er}
H^*(G_\Gamma,\Z)=E/J_K.
\end{equation}

\subsection{Cohomology ring of $\Mo_K$}
\label{subs:macoho}

Recall from Lemma \ref{lem:bg} that, for any compact 
Lie group $G$, there is a fibration sequence
$\Mo_K(EG,G)\to \Mo_K(BG)\to\Mo_K(G^{\times n})$. 
With the additional assumption that $G$ is connected, 
we have an Eilenberg-Moore spectral sequence
\begin{eqnarray}
E_2^{p,q}&=&\Tor^{H^*(BG^{\times n},\Z)}_p(H^*(\Mo_K(BG),\Z),H^*(*,\Z))^q
\nonumber\\
\label{eq:emss}
&=&\Tor^{{H^*(BG^{\times n},\Z)}}_p(H^*(\Mo_K(BG),\Z),\Z)^q
\Rightarrow H^{p+q}(\Mo_K(EG,G),\Z),
\end{eqnarray}
in which resolutions are graded by total degree $q$.  

In the case of $G=S^1$, the cohomology ring $H^*(BG^{\times n},\Z)$ 
is the polynomial ring $S$ with generators in degree $2$.  We have 
a homotopy equivalence $\Mo_K(ES^1,S^1)\simeq\Mo_K(D^2,S^1)$, 
by the discussion in Section~\ref{subsec:mogen}.  As stated in \cite{BP00}, 
the Eilenberg-Moore spectral sequence  \eqref{eq:emss}
degenerates at $E_2$; moreover, one has an isomorphism 
of rings between the cohomology of the moment-angle complex, 
$\Mo_K=\Mo_K(D^2,S^1)$, and the $\Tor$ algebra of $S/I_K$:
\begin{equation}
\label{eq:torkk}
H^*(\Mo_K,\Z)=\Tor^{S} (S/I_K,\Z).
\end{equation}
The argument from \cite{BP00}, completed in \cite{Pa05}, 
relies on the cellular algebra constructed in \cite{BBP04}: 
see the discussion in \S\ref{subsec:cellular}, 
in particular equation \eqref{eq:cwqikoszul}. 

Returning now to the more general situation, let $G$ be an 
arbitrary compact, connected Lie group, and fix a field 
$\k$ of characteristic $0$. By a classical result of Borel, 
$H^*(BG,\k)=H^*(BT,\k)^W$, where $T$ is a maximal torus, 
and $W$ is the Weyl group; consequently, $H^*(BG,\k)$ is 
a polynomial algebra.   It follows that $H^*(BG,\k)$, with $0$ 
as the differential, is a minimal model for $BG$, and so 
$BG$ (and $BG^{\times n}$) are formal spaces; see \cite{Su77} 
and the discussion in \S\ref{subsec:formal}. 

\begin{question}
In this respect, the following interrelated questions arise naturally:
\begin{enumerate}
\item \label{i1}
For which Lie groups $G$ does the Eilenberg-Moore
spectral sequence \eqref{eq:emss} degenerate at $E_2$?
\item \label{i2}
When \eqref{i1} holds, is the
additive isomorphism $E_2\cong H^*(\Mo_K(EG,G),\k)$ an 
isomorphism of rings?
\item \label{i3}
For which groups $G$ is the space $\Mo_K(BG)$ formal?  
(This is established for $G=S^1$ in \cite{NR05}.)
\end{enumerate}
\end{question}

\subsection{Functoriality in $K$}
\label{subsec:funk}

Let $(X,A)$ be a pair of spaces, $L$ a simplicial complex, 
and $K\subset L$ a full subcomplex.  
By Lemma~\ref{lem:full}, the moment-angle complex 
$\Mo_K(X,A)$ is a subcomplex of $\Mo_L(X,A)$, 
with the inclusion map admitting a retraction 
$\Mo_L(X,A)\to \Mo_K(X,A)$.  As a 
consequence, we obtain the following.  

\begin{proposition}
\label{prop:split}
If $K$ is a full subcomplex of $L$, then $H^*(\Mo_K(X,A),\Z)$ 
splits multiplicatively as direct summand of $H^*(\Mo_L(X,A),\Z)$. 
\end{proposition}

One can say more for pairs $(X,*)$  for which $X$ is a 
topological monoid, with $*$ as the unit.    For that, we 
need to switch again to coefficients in our field $\k$.  
By the K\"unneth formula, 
$H^*(X\times X,\k)\cong H^*(X,\k)\otimes H^*(X,\k)$.
It follows that $H^*(X,\k)$ is a Hopf algebra, where 
the multiplication $\mu\colon X\times X\rightarrow X$ 
induces a coproduct $H^*(X,\k)\rightarrow H^*(X,\k)\otimes H^*(X,\k)$ 
given by $x\mapsto x\otimes 1+ 1\otimes x$ for indecomposable $x$; 
see \cite{MM}.

Now suppose that $K$ and $L$ are simplicial complexes on 
$[n]$ and $[m]$, respectively, and $\rho\colon K\rightarrow L$ 
is a map of simplicial complexes.  Then 
$\Mo_\rho\colon \Mo_K(X,X)\rightarrow\Mo_L(X,X)$ 
induces a ring homomorphism 
$\rho^*=H^*(\Mo_\rho)\colon H^*(X,\k)^{\otimes m}\rightarrow 
H^*(X,\k)^{\otimes n}$ 
given by
\begin{equation}
\label{eq:KL}
\rho^*(x^{(j)})=\sum_{i\colon \rho(i)=j}x^{(i)},
\end{equation}
where $x^{(i)}$ denotes the element 
$1\otimes\cdots\otimes1\otimes x\otimes 1\cdots\otimes 1$, 
with $x$ in the $i$th position.  

\begin{lemma}
\label{lem:rhomap}
If $\rho\colon K\rightarrow L$ is a map of simplicial complexes 
as above, and $X=S^1$ or $X=BS^1$, then the homomorphism 
$\rho^*\colon H^*(\Mo_K(X),\k)\rightarrow H^*(\Mo_L(X),\k)$ 
is given by $\rho^*(x_j)=\sum_{i\colon \rho(i)=j}x_i$.  
\end{lemma}

\begin{proof}
From Lemma~\ref{lem:X*quot}, this homomorphism is induced from
$\rho^*\colon H^*(\Mo_K(X,X),\k)\rightarrow H^*(\Mo_L(X,X),\k)$, 
which is given by \eqref{eq:KL}.  The claim follows.
\end{proof}

For $X=BS^1$, this is a result of Panov~\cite{Pa05}.

\section{Homotopy groups and Koszul algebras}
\label{sec:homolie}

In this section, we compute the ranks of the homotopy 
groups of a moment-angle complex $\Mo_K$ in terms 
of the homological algebra of the Stanley-Reisner ring $S/I_K$. 
The answer turns out to be particularly nice in the case when 
$K$ is a flag complex (equivalently, $S/I_K$ is a Koszul algebra). 
For related results, see the work in progress by 
Panov and Ray \cite{PR}. 

Throughout this section, $\k$ continues to denote 
a field of characteristic zero.

\subsection{The homotopy Lie algebra}
\label{subsec:homolie}

For a simply-connected, finite-type CW-complex $X$, 
denote by $\Omega X$ the space of Moore loops at 
the basepoint $*$, endowed with the compact-open 
topology.  The {\em homotopy Lie algebra} of $X$ 
(over $\k$) is the graded $\k$-vector space 
\begin{equation}
\label{eq:homolie}
\g_X= \bigoplus_{r\ge 1}  \pi_r(\Omega X)\otimes \k
\end{equation}
endowed with the graded Lie algebra structure coming 
from the Whitehead product on $\pi_*(X)$ via the 
boundary map in the path fibration $\Omega X\to PX \to X$. 
By the Milnor-Moore theorem \cite{MM}, the universal 
enveloping algebra of $\g_X$ is isomorphic, 
as a Hopf algebra, to $H_*(\Omega X, \k)$.

Now let $K$ be a simplicial complex on vertex set $[n]$ with no 
isolated vertices. From its construction as a subcomplex of 
$BT^n$, it is readily seen that $\Mo_K(BT^n)$ shares the 
same $3$-skeleton with $BT^n=K(\Z^n,2)$;   
hence, $\pi_1(\DJ(K))=0$ and $\pi_2(\DJ(K))=\Z^n$. 
From the long exact homotopy sequence for the fibration 
\begin{equation}
\label{eq:modjbt}
\Mo_K\rightarrow\DJ(K)\rightarrow BT^n, 
\end{equation}
it follows that $\Mo_K$ is $2$-connected, and 
$\pi_q(\Mo_K)\cong\pi_q(\DJ(K))$ for all $q\ge 3$, 
a result of Buchstaber and Panov \cite{BP00}.  
In fact, as these authors note, if $K$ is $k$-neighborly 
(i.e., every $k$-tuple in $[n]$ is a simplex in $K$), 
then $\Mo_K$ is $2k$-connected. 

For convenience, denote the homotopy Lie algebras 
of $\DJ(K)$ and $\Mo_K$ by $\g_{\DJ}$ and $\g_{\Mo}$, 
respectively.  Using the fibration \eqref{eq:modjbt}, we 
obtain a short exact sequence of graded Lie algebras, 
\begin{equation}
\label{eq:omega}
0\to  \g_{\Mo} \to \g_{\DJ} \to L_n\to 0,
\end{equation}
where $L_n$ denotes the abelian Lie algebra of   
rank $n$ generated in degree $1$. Taking Hilbert series 
of enveloping Lie algebras, we obtain
\begin{equation}
\label{eq:hstensor}
h(U(\g_\DJ),t)=h(U(\g_\Mo),t)(1+t)^n.
\end{equation}

\begin{proposition}
\label{prop:ranks}
Let $K$ be a simplicial complex with $n$ vertices 
and $S/I$ the corresponding Stanley-Reisner ring.
Then the ranks $\phi_r=\rank \pi_r(\Mo_K)$ 
of the homotopy groups of the moment-angle 
complex $\Mo_K$ are determined by the following 
identity of formal power series:
\begin{equation}
\label{eq:ranks}
\prod_{r=1}^{\infty} \frac{(1+t^{2r-1})^{\phi_{2r}}}%
{(1-t^{2r})^{\phi_{2r+1}}}=
(1+t)^{-n}\sum_{p,q\geq0}\dim_\k \Tor^{S/I}_{p,q}(\k,\k)\, t^q.
\end{equation}
\end{proposition}

\begin{proof}
By work of Notbohm and Ray \cite{NR05}, the space 
$\DJ(K)$ is known to be formal. Hence, the Eilenberg-Moore 
spectral sequence of the path-loop fibration of $\DJ(K)$ 
degenerates at the $E_2$ term, see \cite{HS79}.  
Together with the Milnor-Moore theorem, this implies 
that the Hopf algebra $U(\g_\DJ)$ is isomorphic to 
$\Tor^{H^*(\DJ(K),\k)}(\k,\k)$.  Finally, since 
$H^*(\DJ(K),\k)=S/I$ (by Theorem~\ref{thm:cohoring}) 
and $\k$ has characteristic zero, we obtain the claimed 
identity from the Poincar\'{e}--Birkhoff--Witt theorem, 
see \cite{MM}. 
\end{proof}

\begin{remark}
\label{rem:golod}
The generating function on the right-hand side of \eqref{eq:ranks} is
known to be rational in the case of a Stanley-Reisner ring: see
\cite{Bac82}.  In the terminology of commutative algebra, the 
ring $S/I$ is {\em Golod-attached,} which means that the canonical  
projection $S\surj S/I$ factors through a sequence of surjections, 
each of which is a {\em Golod homomorphism}.  This in turn 
means that all the Massey products in a certain DGA vanish; 
for a systematic discussion, see Avramov~\cite{Av86b}.
\end{remark}

\subsection{Flag complexes and Koszul duality}
\label{subsec:koszul}

We will give examples later on where 
both sides of \eqref{eq:ranks} can be computed explicitly.  
One such general case is that of $K$ 
a flag complex (meaning, any missing face of $K$ has 
precisely $2$ vertices). 

For a quadratic algebra $A=T(V)/J$, let $A^{!}=T(V^*)/J^{\perp}$ 
denote its quadratic dual.  Let $\Ext_A(\k,\k)$ be the Yoneda 
algebra, and let $\Ext^1_A(\k,\k)$ be the subalgebra generated 
by degree $1$ elements. As shown by L\"{o}fwall \cite{lof86}, 
$A^{!}\cong \Ext^1_A(\k,\k)$.  Recall $A$ is a {\em Koszul algebra} 
if $\Ext^1_A(\k,\k)=\Ext_A(\k,\k)$; see Fr\"oberg's survey \cite{Fr99} 
as a general reference.  Since $\k$ has characteristic zero, 
$\Ext_A(\k,\k)$ is the $\k$-dual of $\Tor^{A}(\k,\k)$.

\begin{theorem}
\label{thm:flag}
Let $K$ be a flag complex, with Stanley-Reisner ring $S/I$.  
Then the ranks $\phi_r$ of the homotopy groups of the 
moment-angle complex $\Mo_K$ are given by:
\begin{equation}
\label{eq:holoflag}
\prod_{r=1}^{\infty} 
\frac{(1+t^{2r-1})^{\phi_{2r}}}{(1-t^{2r})^{\phi_{2r+1}}}=
h(H,i\sqrt{t},-i\sqrt{t})^{-1},
\end{equation}
where $H=\Tor^S(S/I,\k)$ is the cohomology ring of $\Mo_K$, 
and $h(H,s,t)$ is its bigraded Hilbert series.  
\end{theorem}

\begin{proof}
In this case, the ideal $I$ is quadratic, generated by 
monomials $x_ix_j$ for which $\set{i,j}\not\in K$.  Then 
$S/I$ is a Koszul algebra, see \cite{Fr99}.

Now the ``L\"{o}fwall formula'' for quadratic duals says
$h(U(\g_{\DJ}),t)=h(S/I,-t)^{-1}$. On the other hand, 
a standard Euler characteristic calculation shows
\begin{eqnarray}
\label{eq:euler}
h(S/I,t^2) &=& \frac{1}{(1-t^2)^n}\sum_{p,q}(-1)^p{\rm dim}_\k
\Tor^S_{p,q}(S/I,\k) \, t^{p+q}\\
&=&\frac{h(H,t,-t)}{(1-t^2)^n}, \notag
\end{eqnarray}
so $h(U(\g_{\DJ}),t^2)=(1+t^2)^{-n}h(H,it,-it)^{-1}$.
The proof is completed by combining this equality with 
\eqref{eq:hstensor}.
\end{proof}

\section{The cellular cochain algebra and Massey products}
\label{sec:cellmassey}

In this section, we outline Buchstaber and Panov's definition 
of the cochain algebra of $\Mo_K$, and Baskakov's setup 
for computing Massey products in $H^*(\Mo_K,\k)$. We start 
with a discussion of formality and Massey products.  

\subsection{DGA's and formality}
\label{subsec:dga}

Let $\k$ be a field, or the integers. By a {\em differential 
graded algebra}\/ (DGA) we mean a graded  $\k$-algebra $A$, 
endowed with a differential $d\colon A\to A$ of degree $1$.  
The algebra $A$ is said to be {\em (graded)-commutative}\/ 
(CDGA) if $ab =(-1)^{\abs{a}\abs{b}} ba$, for every homogeneous 
elements $a, b\in A$, where $\abs{a}$ denotes the degree of $a$. 

Let $A$ be a DGA. We shall assume the cohomology 
algebra $H^*(A)$ is commutative. We can turn $H^*(A)$ 
into a DGA by assigning to it the zero differential.  
The differential graded algebra $(A,d)$ is said to be 
{\em formal} if there is a sequence of DGA morphisms  
(going either way), connecting $(A,d)$ to $(H^*(A),0)$ 
and inducing isomorphisms in cohomology.  

The simplest formality test is provided by the Massey 
products.  As is well-known, if $(A,d)$ is formal, then 
all Massey products (of order $3$ or higher) vanish; 
see \cite{DGMS, Su77}. 
Let us briefly review the relevant definitions; for 
simplicity, we will only treat the triple products here. 

Assume $\alpha_1,\alpha_2,\alpha_3$ 
are homogeneous elements in $H^*(A)$ such 
that $\alpha_1 \alpha_2 =\alpha_2\alpha_3=0$.  
Pick representative cocycles $a_i$ for $\alpha_i$, 
and elements $x,y\in A$ such that $dx=a_1a_2$ 
and $dy=a_2a_3$.  Setting $\bar{a}=(-1)^{\abs{a}}$, it 
is readily seen that  $xa_3-\bar{a}_1y$ is a cocyle.  
The set of cohomology classes of all 
such cocycles is the Massey triple product
$\angl{\alpha_1,\alpha_2,\alpha_3}$.   
The image of this set  in the quotient ring 
$H^*(A)/( \alpha_1,\alpha_3)$ is a well-defined 
element of degree 
$|\alpha_1|+|\alpha_2|+|\alpha_3|-1$; we say 
$\angl{\alpha_1,\alpha_2,\alpha_3}$ 
is {\em non-vanishing} if this element is not $0$.  
The Massey product $\angl{\alpha_1,\alpha_2,\alpha_3}$ 
is said to be {\em decomposable} if it contains a cohomology 
class that can be written as a product $\lambda \nu$ 
of two elements in $H^{>0}(A)$; otherwise, it is called 
indecomposable.  

\begin{lemma}
\label{lem:cellmas}
Let $(A,d)$  and $(A',d')$ be two DGA's.  Suppose 
$\rho\colon (A',d') \to (A,d)$ is a surjective chain map, 
chain homotopic to a ring map, and inducing an \nopagebreak
isomorphism $\rho^*\colon H^*(A',d') \to H^*(A,d)$.  
If $\angl{\alpha_1,\alpha_2,\alpha_3}$ 
is a non-vanishing Massey product in $H^*(A,d)$, then 
$(\rho^*)^{-1} \angl{\alpha_1,\alpha_2,\alpha_3}$ is 
 a non-vanishing Massey product in $H^*(A',d')$.  
\end{lemma}

\begin{proof}
Pick representatives as above, so that 
$\angl{\alpha_1,\alpha_2,\alpha_3} =[xa_3-\bar{a}_1y]$.  
Using the fact that $\rho$ is a surjective quasi-isomorphism, 
we may find cocycles $a'_1,a'_2,a_3'\in  A'$ such that 
$\rho(a'_i)=a_i$, and cochains $x',y'\in A'$ such that 
$dx'=a'_1a'_2$ and $dy'=a'_2a'_3$, and thus form the 
Massey product 
$\angl{\alpha'_1,\alpha'_2,\alpha'_3} =[x'a'_3-\bar{a}'_1y']$.  

Now, by assumption, there is a a degree $1$ map 
$\sigma\colon (A',d') \to (A,d)$ so that 
$\rho(a'b')-\rho(a')\rho(b')=(\sigma d'+d \sigma)(a'b')$, 
for every $a',b'\in A'$.   This implies 
$\rho(x'a'_3-\bar{a}'_1y')-(xa_3-\bar{a}_1y) \in \im(d)$, 
finishing the proof.
\end{proof}

\subsection{Formal spaces}
\label{subsec:formal}

To a space $X$, Sullivan associates in \cite{Su77} 
a commutative DGA: the rational algebra of polynomial 
forms, $(A_{\rm PL}(X;\Q),d)$.  The space $X$ is said 
to be {\em rationally formal}\/ (or simply, formal) if 
$(A_{\rm PL}(X;\Q),d)$ is formal in the category 
of CDGA's.   For more details, and equivalent 
definitions, see \cite{DGMS, FHT}. 

Examples of formal spaces include spheres; Eilenberg-MacLane 
spaces $K(\pi,n)$ with $n>1$; compact, connected Lie groups $G$ 
and their classifying spaces $BG$ \cite{Su77};  and compact K\"{a}hler 
manifolds \cite{DGMS}.  Formality is preserved under wedges and 
products of spaces, and connected sums of manifolds.

We say $X$ is {\em integrally formal}\/ if the singular 
cochain algebra $(C^{*} (X,\Z),d)$ is formal.  
From the definition it is apparent that, if $(C^*(X,\Z),d)$ is
formal, so is $(C^*(X,\k),d)$, for any choice of field $\k$ 
of characteristic $0$.  A result of Watkiss, recorded 
in \cite[Corollary 10.10]{FHT}, implies that if $X$ is 
rationally formal, then $(C^*(X,\k),d)$ is formal, 
where $\k$ is any field of characteristic zero.

In what follows,  we will show that, for certain spaces $X$,
the singular cochain DGA $(C^*(X,\Q),d)$ is not formal.  By the
logic above, then, $X$ is neither integrally formal, nor 
rationally formal.  Since the notions coincide, we will 
simply write that $X$ is not formal.

In \cite{NR05}, Notbohm and Ray show that the 
Davis-Januszkiewicz spaces, $\DJ(K)=\Mo_K(BS^1)$, 
are both integrally formal and rationally formal. 
It is also known that  the classifying spaces 
for right-angled Artin groups, 
$K(G_{\Gamma},1)=\Mo_{\Delta(\Gamma)} (S^1)$, 
are rationally formal; see \cite{PS06}.

If $X$ is formal, and $\angl{\alpha_1,\alpha_2,\alpha_3}$ 
is a Massey product computed from $(C^*(X,\k),d)$, 
then $\angl{\alpha_1,\alpha_2,\alpha_3}$ vanishes.   
Thus, non-vanishing Massey products provide a 
handy tool for detecting non-formality.  

As shown by Stasheff \cite{Sta83}, if $X$ is a $k$-connected 
CW-complex of dimension $n \le 3k+1$, then $X$ is formal.  
This is best possible:  attaching a cell $e^{3k+2}$ to the 
wedge $S^{k+1}\vee S^{k+1}$ via the iterated Whitehead 
product $[\iota_1,[\iota_1,\iota_2]]$ yields a non-formal 
CW-complex.  According to Miller \cite{Mi79}, the dimension 
bound can be relaxed for manifolds:  if $X$ is a compact, 
$k$-connected manifold  of dimension $n \le 4k+2$, 
then $X$ is formal.   Again, this bound is best possible, 
see \cite{FM04}, \cite{DR05}.

\subsection{The cellular cochain algebra of $\Mo_K$}
\label{subsec:cellular}

Given a CW-complex $X$, let $(C^*_{\CW}(X),d)$ be its 
(integral) cellular cochain complex.  The cellular cup-product 
map may be defined as the composite 
\begin{equation}
\label{eq:cellcup}
\xymatrix{
\cup \colon\, C^*_{\CW}(X)\otimes C^*_{\CW}(X) \ar^(.6){\times} [r]  & 
C^*_{\CW}(X\times X) \ar^(.55){\widetilde\Delta} [r]  &
C^*_{\CW}(X),
}
\end{equation}
where $\widetilde\Delta$ is a cellular approximation to the 
diagonal map $\Delta \colon X\to X\times X$.   Unfortunately, 
such an approximation is not functorial, and the cup-product 
map on $C^*_{\CW}(X)$ is not associative, in general---although, 
of course, $\cup$ induces an associative and commutative 
product on $H^*(X,\Z)$.

Despite these difficulties, Buchstaber, Panov, and Baskakov 
\cite{BP00, BBP04, Pa05} were able to define a DGA structure 
on the cellular cochain complex of a moment-angle 
complex $\Mo_K=\Mo_K(D^2,S^1)$, over a field $\k$. 
Let us review their construction.  

The unit disk $D^2=\{z\in \C \mid \abs{z}\le 1\}$ comes 
endowed with a natural cell structure:  $D^2=e^0\cup e^1 \cup e^2$, 
with $0$-skeleton the point $z=1$, with $1$-skeleton the boundary 
circle $S^1$, and with $2$-skeleton the disk itself.  The 
diagonal map of $D^2$ has a natural cellular approximation, 
$\widetilde\Delta \colon D^2 \to D^2\times D^2$, sending 
$z=re^{i\theta}$ to  
$(1+r(e^{2i\theta}-1),1)$  for $0\le \theta< \pi$ and to 
$(1,1+r(e^{2i\theta}-1))$  for $\pi\le \theta< 2\pi$. 
With the cup product defined by $\widetilde\Delta$, 
the cellular cochain complex $C^*_{\CW}(D^2)$ 
becomes a commutative DGA, with generators $u$ and $x$ 
dual to $e^1$ and $e^2$, multiplication $x^2=xu=0$, and 
differential $du=x$, $dx=0$.   

Now let $K$ be a simplicial complex on $[n]$, and let $\Mo_K\subset 
(D^2)^{\times n}$ be the corresponding moment-angle complex.  
The diagonal map of $(D^2)^{\times n}$ has cellular approximation 
$\widetilde\Delta^{\times n}$; this map restricts to a cellular 
approximation $\widetilde\Delta_K$ to the diagonal map of $\Mo_K$.   
The resulting cup product on the cellular cochain complex of 
$\Mo_K$, with coefficients in $\k$, yields a commutative DGA 
with presentation 
\begin{equation}
\label{eq:cellzk}
C^*_{\CW}(\Mo_K,\k) = 
\big( (S/I_K) \otimes_{\k} E \big) / (x_i^2=x_iu_i=0),
\end{equation}
where $S=\k[x_1,\dots, x_n]$ is the polynomial ring with 
generators $x_i$ in degree $(0,2)$, 
$E=\bigwedge [x_1,\dots, x_n]$ is the exterior algebra 
with generators $u_i$ in degree $(1,1)$, 
and with differential $d$ given by $du_i=x_i$, $dx_i=0$.  
Additively, this algebra
is generated by square-free monomials $x_\sigma u_{\I}$, 
for all choices of $\sigma\in K$ and $\I\subseteq[n]$ for 
which $\sigma\cap \I=\emptyset$.
The construction has the following naturality property:  
If $f\colon K\to L$ is a simplicial map, then 
$\Mo_f\colon \Mo_K\to \Mo_L$ commutes 
with the diagonal approximations, and so 
$\Mo_f^* \colon C^*_{\CW}(\Mo_L,\k)  \to C^*_{\CW}(\Mo_K,\k)$  
is a DGA morphism.  

Baskakov, Buchstaber and Panov~\cite{BBP04} show that the natural surjection 
from the Koszul complex of the Stanley-Reisner ring to the 
cellular cochain algebra, 
\begin{equation}
\label{eq:cwqikoszul}
\pi\colon S/I_K\otimes_{\k} E \rightarrow
C^*_{\CW}(\Mo_K,\k),
\end{equation}
is a (bigraded) quasi-isomorphism of DGAs.  
This recovers the isomorphism 
$H^q(\Mo_K,\k)_p\cong \Tor_{p,q}^S(S/I_K,\k)$ from \eqref{eq:torkk}, 
grading $\Tor$ by homological and internal degree, respectively.

\subsection{Multiplicative structures in resolutions}
\label{subsec:multres}

Formality of the moment-angle complex $\Mo_K$ has 
a quite closely related notion in commutative algebra.  

Given a $\k$-algebra homomorphism 
$S\to S/I$, a projective resolution $P$ of $S/I$ 
over $S$ is said to have {\em multiplicative structure} if $P$ 
has the structure of a commutative DGA, and the surjection
$P\surj S/I$ is a quasi-isomorphism of DGAs over $S$   
(for details, see Avramov's survey \cite{Av98}).  
Also recall that graded modules over commutative 
graded rings have minimal free resolutions:  that is, 
resolutions in which the differential can be expressed as a 
matrix whose entries have strictly positive degree.  
Such resolutions are essentially unique.

\begin{proposition}
\label{prop:minres}
Suppose the minimal free resolution $P$ of $S/I$ over $S$ has
multiplicative structure.  Then $C^*_{\CW}(\Mo_K,\k)$ is formal.
\end{proposition}

\begin{proof}
The Koszul complex $S\otimes_\k E$ is a DGA resolution 
of $\k$. Since $\Tor^S(S/I,\k)$ may be computed by 
resolving either $S/I$ or $\k$, the augmentation maps 
in each resolution give quasi-isomorphisms of DGAs
(see \cite[2.3.2]{Av98}):
\begin{equation}
\label{eq:minresqi}
P\otimes_S\k\leftarrow P\otimes_S E\rightarrow S/I\otimes_S E,
\end{equation}
each having homology $\Tor^S(S/I,\k)=H^*(\Mo_K,\k)$. 
Since $P$ is minimal, however, $P\otimes_S\k$ also has trivial 
differential, so $P\otimes_S\k\cong H^*(\Mo_K,\k)$, which 
means that each DGA in \eqref{eq:minresqi} is formal.  
By \eqref{eq:cwqikoszul}, $C^*_{\CW}(\Mo_K,\k)$ is formal.
\end{proof}

\subsection{The Baskakov formula}
\label{subsec:bask}

Using the cellular cochain algebra model for $\Mo_K$, 
Baska\-kov \cite{Ba02} gave an explicit formula for the 
cup product in  $H^{*}(\Mo_K,\Z)$, in terms of pairings 
between full subcomplexes.  Let us recall Baskakov's 
formula (see also Panov \cite{Pa05}) by means of \eqref{eq:cellzk}.

For each subset $\I\subseteq[n]$, let $K_\I$ denote the full  
subcomplex of $K$ on vertices $\I$. For each $\sigma\in K_\I$, 
let $s=\abs{\sigma}$ and $p=\abs{\I}$.  Then define a map to 
the (reduced) group of simplicial cochains on $K_\I$,
\begin{equation}
\label{eq:phicw}
\phi_\I\colon C^{p,p+2s}_{\CW}(\Mo_K)\rightarrow 
\widetilde{C}^{s-1}(K_\I),
\end{equation}
by letting $\phi_\I(x_\sigma u_{\I-\sigma})=\chi_\sigma$,
where $\chi_\sigma$ denotes the indicator function on 
the simplex $\sigma$ in $K_\I$. (By convention, 
$\widetilde{C}^{-1}(\emptyset)=\k$.) 
The maps $\set{\phi_\I}$ assemble to give an isomorphism 
of cochain complexes, 
\begin{equation}
\label{eq:Phimap}
\Phi \colon C^*_{\CW}(\Mo_K)\rightarrow\bigoplus_{\I\subseteq[n]}
\widetilde{C}^*(K_\I).
\end{equation}

The map $\Phi$ is, in fact, a ring isomorphism, 
where the multiplication on the left is defined by  
\eqref{eq:cellzk}, while the multiplication on the right 
is defined as follows. 
For $\alpha\in \widetilde{C}^{p}(K_{\I})$ and $\alpha'\in 
\widetilde{C}^{p'}(K_{\I'})$, the product $\alpha\cdot\alpha'\in
\widetilde{C}^{p+p'+1}(K_{\I\sqcup \I'})$ is given by
\begin{equation}
\label{eq:baskakov}
\alpha\cdot\alpha'=
\begin{cases}
0, &
\text{if  $\I\cap \I'\ne\emptyset$} ,\\[3pt]
j^* \psi (\alpha \otimes \alpha'), &
\text{if  $\I\cap \I'=\emptyset$} ,
\end{cases}
\end{equation}
where 
$ j \colon K_{\I\sqcup \I'}=K_{\I}\sqcup
K_{\I'}\hookrightarrow K_{\I}*K_{\I'}$ 
is the inclusion into the simplicial join, and 
$\psi \colon\widetilde{C}^{p}(K_{\I})\otimes
\widetilde{C}^{p'}(K_{\I'})  \xrightarrow{\:\cong\:} 
\widetilde{C}^{p+p'+1}(K_{\I}*K_{\I'})$ 
is the standard isomorphism between the respective 
reduced simplicial cochain complexes.

Taking cohomology, one obtains the (additive) formula 
of Hochster~\cite{Ho77},
\begin{equation}
\label{eq:hochster}
\Tor^S_{p-s,p+s}(S/I,\Z) = 
\bigoplus_{\stackrel{I\subseteq[n]}{\abs{I}=p}}
\widetilde{H}^{s-1}(K_\I,\Z).
\end{equation}

\subsection{Massey products in $\Mo_K$}
\label{subsec:zkm}

Let us now show that the non-formality of $\Mo_K$ can be 
detected by non-vanishing Massey products---not just in the 
singular cochain algebra $C^{*}(\Mo_K,\k)$, but also in the 
much simpler (and commutative) cellular cochain 
algebra $C^{*}_{\CW}(\Mo_K,\k)$.  

\begin{proposition}
\label{prop:cellmas}
Suppose the differential graded algebra $C^*_{\CW}(\Mo_K,\k)$ 
carries a non-vanishing triple Massey product.
Then $\Mo_K$ is not formal.
\end{proposition}

\begin{proof}
Simplicial approximation yields a surjective homomorphism 
$\rho\colon C^*((D^2)^{\times n},\k) \surj 
C^*_{\CW}((D^2)^{\times n},\k)$, commuting with the 
differentials.  Since $(D^2)^{\times n}$ 
is acyclic,  $\rho(ab)-\rho(a)\rho(b)$ is a cellular coboundary, 
for any singular cochains $a,b$; thus, $\rho$ is also 
a ring map, up to chain homotopy.  Now consider the 
commuting diagram
\begin{equation}
\label{eq:cdsing}
\xymatrix{
C^*((D^2)^{\times n},\k)  \ar@{->>}^{\rho\:}[r] \ar@{->>}^{\pi}[d] 
&C^*_{\CW}((D^2)^{\times n},\k)\ar@{->>}^{\pi}[d] \\
C^*(\Mo_K,\k) \ar@{->>}^{\rho_K}[r] &C^*_{\CW}(\Mo_K,\k) 
}
\end{equation}
where the vertical arrows are the canonical projections 
(both morphisms of DGA's), while $\rho_K$ is the restriction 
of $\rho$.  Chasing the diagram, we see that the chain map 
$\rho_K$ is also a ring map, up to chain homotopy.  The 
conclusion follows from Lemma \ref{lem:cellmas} and the 
standard Massey product formality test for $C^{*}(\Mo_K,\k)$. 
\end{proof}

\section{Triple Massey products in lowest degree}  
\label{sec:decomp}

In \cite{Ba03}, Baskakov constructs a family of simplicial 
complexes $K$ for which the cellular cochain algebra 
defined above has non-vanishing Massey triple products: 
a complete discussion is found in \cite{Pa05}.  Here, we 
show that one can characterize the complexes $K$ for 
which $H^*(\Mo_K,\k)$ contains a non-trivial Massey product 
in lowest possible degree.  

\subsection{The five obstruction graphs}
\label{subsec:five obs}
Let $K$ be a simplicial complex on $n$ vertices.  Let 
$\ones$ denote the $1$-skeleton of $K$.

\begin{theorem}
\label{th:char}
The following are equivalent:
\begin{romenum}
\item 
There exist cohomology classes 
$\alpha,\beta,\gamma\in H^3(\Mo_K,\k)$ for which
$\angl{\alpha,\beta,\gamma}$ is defined and non-trivial. 
\item 
The underlying graph of $\ones$ contains an induced 
subgraph isomorphic to one of the five graphs in 
Figure~\ref{fig:exclude}.  
\end{romenum}
Moreover, all Massey products arising in this fashion are 
decomposable.
\end{theorem}

\begin{figure}
\setlength{\unitlength}{1cm}
\begin{picture}(8,2.5)(1,-2)
\xygraph{
\no="4"-[]!{!(0,-0.5)}\no="6"-[]!{!(0,-0.5)}\no="3"
-[]!{!(0.5,0)}\no="5"-[]!{!(0,0.5)}\no="2"-"6"-[]!{!(-0.65,-0.9)}\no="1"
-@/_/"5"
"1"-@/^/"4"-"2"
"1"-"3"
}
\quad 
\xygraph{
\no="4"-[]!{!(0,-0.5)}\no="6"-[]!{!(0,-0.5)}\no="3"
-[]!{!(0.5,0)}\no="5"-[]!{!(0,0.5)}\no="2"-"6"[]!{!(-0.65,-0.9)}\no="1"
-@/_/"5"
"1"-@/^/"4"-"2"
"1"-"3"
}
\quad
\xygraph{
\no="4"-[]!{!(0,-0.5)}\no="6"-[]!{!(0,-0.5)}\no="3"
-[]!{!(0.5,0)}\no="5"-[]!{!(0,0.5)}\no="2""6"-[]!{!(-0.65,-0.9)}\no="1"
-@/_/"5"
"1"-@/^/"4"-"2"
"1"-"3"
}
\quad
\xygraph{
\no="4"-[]!{!(0,-0.5)}\no="6"-[]!{!(0,-0.5)}\no="3"
-[]!{!(0.5,0)}\no="5"-[]!{!(0,0.5)}\no="2""6"[]!{!(-0.65,-0.9)}\no="1"
-@/_/"5"
"1"-@/^/"4"-"2"
"1"-"3"
}
\quad
\xygraph{
\no="4"-[]!{!(0,-0.5)}\no="6"-[]!{!(0,-0.5)}\no="3"
-[]!{!(0.5,0)}\no="5"-[]!{!(0,0.5)}\no="2""6"[]!{!(-0.65,-0.9)}\no="1"
"1"-@/^/"4"-"2"
"1"-"3"
}
\end{picture}
\caption{\textsf{The five obstruction graphs}}
\label{fig:exclude}
\end{figure}

By the functoriality of the moment-angle construction on 
inclusions of full subcomplexes $L\hookrightarrow K$, 
(Lemma~\ref{lem:full} and subsection~\ref{subsec:cellular}), 
the Theorem is implied by the following statement:

\begin{proposition}
\label{prop:char}
If $L$ is a simplicial complex on six vertices, $H^*(\Mo_L,\k)$ has a 
non-trivial Massey triple product if and only if its $1$-skeleton is
shown in Figure~\ref{fig:exclude}.  Moreover, any non-trivial Massey 
product in $H^*(\Mo_L,\k)$ is decomposable.
\end{proposition}

The proof of this Proposition will occupy the rest of this section.

\subsection{Proof of ``$\Leftarrow$'' of Proposition~\ref{prop:char}:}
\label{subsec:necessary}
Suppose the one-skeleton of $L$ is one of Figure~\ref{fig:exclude}.
None of these graphs contains a complete subgraph on four vertices,
from which it follows $L$ is a $2$-complex.  The only two-cells which 
may appear are the interiors of the triangles.  Regardless of the 
two-cells, $\widetilde{H}^p(L,\k)\neq 0$ iff $p=1$.

Label the vertices of the ``obstruction'' graphs as in Figure 
\ref{fig:excl}. Then $L^{(1)}$ is obtained by omitting any 
subset of the dotted edges.  Recall $L_{\I}$ denotes the full 
subcomplex on index set $\I$; we will slightly abuse notation, 
and write, e.g., $L_{i,j}$ or even $L_{ij}$ for $L_{\{i,j\}}$.  

\begin{figure}[ht]
\setlength{\unitlength}{1cm}
\begin{picture}(4,4)(-0.5,-3.6)
\xygraph{
*+{4}="4"-[d]*+{6}="6"-[d]*+{3}="3"
-[r]*+{5}="5"-[u]*+{2}="2"-@{.}"6"-@{.}[]!{!(-1.3,-1.8)}*+{1}="1"
-@/_/@{.}"5"
"1"-@/^/"4"-"2"
"1"-"3"
}
\end{picture}
\caption{\textsf{A labeled obstruction graph}}
\label{fig:excl}
\end{figure}

Let $\alpha$, $\beta$, and $\gamma$ be non-trivial classes in 
$\widetilde{H}^0(L_{i,i+1},\k)$, for $i=1,3,5$, respectively.  Abusing
notation, we shall identify $\alpha$, $\beta$, and $\gamma$ with
their images in $H^3(\Mo_L,\k)$ under the isomorphism $\Phi_*$ 
induced by the chain map \eqref{eq:Phimap}.  Then, by means of 
Baskakov's formula \eqref{eq:baskakov}, we show the triple product 
$\angl{\alpha,\beta,\gamma}\in H^8(\Mo_L,\k)$ is non-trivial, as follows.

We will work with cellular cochains and denote cocycles by indicator
functions on simplicial subcomplexes.  It should be noted that the 
underlying subcomplex is part of the data, but for legibility reasons 
will be taken to be implicit in what follows.

Since Massey products are linear in each argument, we may replace
each cocycle with a nonzero scalar multiple without loss of generality.
Accordingly, let $\alpha=[\chi_2]$, $\beta=[\chi_3]$, and 
$\gamma=[\chi_5]$.  From Figure \ref{fig:excl}, we see the 
subcomplexes $L_{1234}$ and $L_{3456}$ are both paths: 
in particular 
\begin{equation}
\label{eq:contr}
\widetilde{H}^p(L_{1234},\k)=
\widetilde{H}^p(L_{3456},\k)=0\quad\hbox{for all $p$},
\end{equation}
so $\alpha\beta=\beta\gamma=0$.

Then, multiplying cochains with the Baskakov construction, 
we see $\chi_2\chi_3=0$ while
$\chi_3\chi_5=\chi_{35}=-\delta(\chi_5)$. 
It follows that the triple product is represented by the cochain
\begin{equation}
\label{eq:chi25}
0\cdot \chi_5-(-1)^3\chi_2\cdot(-\chi_5) = -\chi_{25}.
\end{equation}
Let $\omega=[-\chi_{25}]$.  It is easily seen that $\omega$ is 
nonzero in $H^1(L_{123456},\k)=H^8(\Mo_L,\k)$, 
regardless of the filling of the two-cells.

Last, we must check that the Massey product 
$\angl{\alpha,\beta, \gamma}$ is non-trivial with respect 
to indeterminacy; that is, $\omega$
does not lie in the ideal of $H^*(\Mo_L,\k)$ generated by 
$\set{\alpha,\gamma}$.  For this, suppose $\alpha\eta\in H^8(\Mo_L,\k)$ 
is supported on $L_{123456}$.  For grading reasons, the support of 
$\eta$ is $L_{3456}$.  However, from \eqref{eq:contr}, $\eta$ must 
be zero.  The same argument applies to $\gamma$, so 
$(\alpha,\gamma)\cap H^8(\Mo_L,\k)=0$;
the indeterminacy for this Massey product is zero.

Finally, we claim that $\omega=\nu\nu'$ for certain 
classes $\nu,\nu' \in H^4(\Mo_L,\k)$; 
that is, $\omega$ is decomposable.  Let $\nu=[-\chi_2]$
on $L_{123}$, a generator of $\widetilde{H}^0(L_{123},\k)$.  
Similarly, let $\nu'=[\chi_5]$ on $L_{456}$.  Since 
$\omega=[-\chi_{25}]$, the claim follows by \eqref{eq:baskakov}.

\subsection{Proof of ``$\Rightarrow$'' of Proposition~\ref{prop:char}:}
\label{subsec:sufficient}
Suppose that $L$ is a simplicial complex on six vertices, and 
it possesses a non-trivial Massey product $\angl{\alpha,\beta,\gamma}$.  
For grading reasons, each of $\alpha,\beta,\gamma$ must be supported on 
pairwise-disjoint, $2$-vertex, disconnected subcomplexes.  Without
loss of generality, suppose $\alpha,\beta,\gamma$ are supported on
vertices $S_1=\set{1,2}$, $S_2=\set{3,4}$, and $S_3=\set{5,6}$, 
respectively.  

Let the graph $G$ be the edge-complement of $L^{(1)}$.  
We shall show that $G$ must be obtained from the graph 
below by adding any subset of the dotted edges.  Then, 
by comparing with Figure \ref{fig:excl}, the implication
is proven.
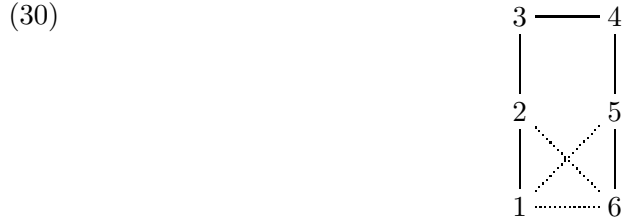
\begin{equation}
\label{eq:complgraph2}
\xygraph{
*+{3}="3" -[r] *+{4} -[d] *+{5}="5" -[d] *+{6}="6"
-@{.}[l] *+{1}="1" -[u] *+{2}="2" 
"5"-@{.}"1"
"2"-@{.}"6"
"3"-"2"
}
\end{equation}

Immediately, $G$ contains edges $i,i+1$ for $i=1,3,5$. 

\begin{lemma}
\label{lem:g12}
$G$ contains edges joining $S_i$ to $S_{i+1}$ for $i=1,2$.
\end{lemma}

\begin{proof}
Since $\alpha,\beta,\gamma$ define a Massey product, $\alpha\beta=0$.
It follows that $L_{1234}$ is not a cycle, so there is at least one
edge from $S_1$ to $S_2$.  A symmetric argument applies to $S_2,S_3$.
\end{proof}

\begin{lemma}
\label{lem:g123}
$G$ does not contain a path $v_1v_2v_3$, where $v_i\in S_i$ for $i=1,2,3$.
\end{lemma}

\begin{proof}
Suppose it did.  Then, up to nonzero scalar multiples,
$\alpha=[\chi_{v_1}]$, $\beta=[\chi_{v_2}]$, and $\gamma=[\chi_{v_3}]$. 
The edge $v_iv_{i+1}$ in $G$ means there is no edge between $v_i$ and $v_{i+1}$
in $L$, so $\chi_{v_i}\chi_{v_{i+1}}=0$ for $i=1,2$.  Then
the Massey product vanishes, contradicting the hypotheses.
\end{proof}

Up to relabeling, then, $G$ contains the path $12\cdots6$, and vertices
$3,4$ both have degree $2$.  It remains to check the following.

\begin{lemma}
\label{lem:g25}
$G$ does not contain the edge $25$.
\end{lemma}

\begin{proof}
If it did, represent $\alpha,\beta,\gamma$ by $\chi_2$, $\chi_3$, $\chi_5$,
respectively.  Then $\chi_2\chi_3=0$, while $\chi_3\chi_5=-\delta(\chi_5)$.
Then an edge $v_2v_5$ would make $\chi_2\chi_5=0$ and the Massey 
product cohomologous to zero, a contradiction.
\end{proof}

The proof of Proposition~\ref{prop:char} now follows directly from the
previous three Lemmas.  The Massey product $\angl{\alpha,\beta,\gamma}$ 
is, up to a nonzero scalar, the one constructed in \ref{subsec:sufficient};
in particular, it is decomposable.

\section{Formal moment-angle manifolds}
\label{sec:formalman}

In this and the next section, we consider closed manifolds obtained 
from sphere triangulations by means of the moment-angle 
construction.  We begin with the formal ones.

\subsection{Moment-angle manifolds}
\label{subsec:momman}

Suppose $K$ is an $n$-vertex triangulation of the sphere 
$S^{\ell}$.  Then, as shown by Buchstaber and Panov \cite{BP00}, 
the moment-angle complex $\Mo_K=\Mo_K(D^2,S^1)$ is 
a smooth, compact, connected manifold of dimension $n+\ell+1$; 
moreover, $\Mo_K$ is $2$-connected, and in fact, $2k$-connected, 
if $K$ is $k$-neighborly.

Some moment-angle manifolds can be described 
in simple terms.  As noted in Example \ref{ex:mom}, if 
$K$ is the boundary of the $(n-1)$-simplex, then  
$\Mo_K=S^{2n-1}$.   Also, if $K=(S^0)^{*d}$ 
(the simplicial join of $d$ copies of $S^{0}$) 
is the $(d-1)$-dimensional hyperoctahedron, 
then $\Mo_K=(S^3)^{\times d}$, from Corollary~\ref{lem:joins}.  
In general, though, moment-angle manifolds can exhibit quite 
a complicated structure, both from the point of view of 
their cohomology ring and their Massey products. 

We start with an observation concerning the formality of 
such manifolds.

\begin{proposition}
\label{prop:formalman}
Suppose $K$ is an $n$-vertex triangulation of $S^{\ell}$. 
\begin{romenum}
\item If $n+\ell \le 9$, then $\Mo_K$ is formal.  
\item If $K$ is $k$-neighborly, and $n+\ell\le 8k+1$, then 
 $\Mo_K$ is formal. 
 \end{romenum}
\end{proposition}

\begin{proof}
Use  Miller's result \cite{Mi79}, as recounted 
in \S\ref{subsec:formal}.   
\end{proof} 

In particular, all triangulations of $S^2$ on at most 
$7$ vertices (there are precisely $9$ such) 
yield formal  moment-angle manifolds.

\subsection{Corner-cutting}
\label{subsec:corner}

We now describe an operation on simplicial complexes which 
proves to be useful in this context.     

Given a simplicial complex $K$ and a maximal simplex 
$F=(v_0,v_1,\ldots,v_k)$ in $K$, let $K*_F w$ denote 
the complex obtained by adding a new vertex $w$, 
removing the simplex $F$, and adding in $k+1$ 
maximal simplices $(v_0,v_1,\ldots,v_{i-1},w,v_{i+1},\ldots,v_k)$ 
for $0\leq i\leq k$.

Note that this operation is dual to cutting the corner 
from a polytope:   if $K$ triangulates $S^k$ and 
$F$ is a maximal simplex of $K$, then let $P$ 
and $P'$ denote the dual polytopes to $K$ and 
$K*_F w$, respectively.  Maximal simplices of 
the triangulation are in bijection with vertices 
of the dual polytope.  So $P'$ is obtained from 
$P$ by cutting off the vertex labeled by $F$ in $P$.  

Work of McGavran \cite[Theorem 3.4]{McG79}  and 
Bosio and Meersseman \cite[Theorem~6.3]{BM04} 
gives the following:  if $K$ is obtained from the boundary 
of a simplex by a sequence of moves as above, then 
$\Mo_K$ is diffeomorphic to 
\begin{equation}
\label{eq:csum}
\cs_{i=1}^{p} S^{a_i}\times S^{n-a_i},
\end{equation}
a connected sum of products of spheres.  Consequently, 
such moment-angle manifolds are formal spaces, and 
the cohomology ring $H^*(\Mo_K,\Z)$ has a very simple 
structure, essentially dictated by the Betti numbers 
and Poincar\'e duality.  

\begin{example}
\label{ex:polygons}
Denote by $\Mo_n$ the $(n+2)$-dimensional moment-angle 
manifold corresponding to an $n$-gon, $n\ge 3$.   
From the discussion in \S\ref{subsec:mogen}, we know that 
$\Mo_3=S^5$ and $\Mo_4=S^3\times S^3$.  More generally, 
by \cite{McG79} we have:
\begin{equation}
\label{eq:polycsum}
\Mo_n = \cs_{j=1}^{n-3} \,  j \tbinom{n-2}{j+1} S^{j+2} \times S^{n-j} ,
\end{equation}
for all $n\ge 4$.  Computing Betti numbers from this 
decomposition, we find  
\begin{equation}
\label{eq:bettizn}
b_k(\Mo_n)= \frac{k(k-2)(n-k)}{(n-1)(n-k+1)} \cdot \binom{n}{k} ,
\end{equation}
for $3\le k\le n-1$, which recovers a result of Buchstaber 
and Panov \cite[Example 4.3.5]{BP00}.

Now, for all $n\ge 4$, the $n$-gon is a flag complex, 
so Theorem \ref{thm:flag} applies to describe the ranks of the
homotopy groups of $\Mo_n$.  We will do so directly.
Using a result of Stanley~\cite{St96}, the Hilbert series
of $S/I$ is given by the number of faces of each 
dimension in the $n$-gon:
\begin{equation}
\label{eq:hilbzn}
h(S/I,t)=1+\frac{nt}{1-t}+\frac{nt^2}{(1-t)^2}.
\end{equation}
Then, by \eqref{eq:euler}, 
\begin{eqnarray}
\label{eq:htt}
h(H,t,-t)&=&(1-t^2)^nh(S/I,t^2)\\
&=&(1-t^2)^{n-2}(1+(n-2)t^2+t^4),\notag
\end{eqnarray}
where we recall 
$h(H,s,t)=\sum_{p,q}\dim_\k \Tor^S_{p,q}(S/I,\k) \, t^ps^q$.
Using \eqref{eq:holoflag}, 
\begin{equation}
\label{eq:phin}
\prod_{r=1}^{\infty} 
\frac{(1-t^{2r})^{\phi_{2r+1}}}{(1+t^{2r+1})^{\phi_{2r+2}}}=
(1+t)^{n-2}(1+(2-n)t+t^2),
\end{equation}
from which the ranks of the homotopy groups may be computed recursively.
For example,
\begin{eqnarray}
\label{eq:phirzn} 
\phi_3 &=&n(n-3)/2 \notag \\
\phi_4 &=&n(n-2)(n-4)/3 \notag \\
\phi_5 &=&n(n-1)(n-3)(n-4)/4\\ 
\phi_6 &=&n(n-1)(n-2)(n-3)(n-4)/5 \notag \\
\phi_7 &=&n(n-2)(n-3)(n-4)(n^2-3n+1)/6, \notag
\end{eqnarray}
for $n\geq4$.  
\end{example}

\subsection{Beyond connected sums of products of spheres}
\label{subsec:notcs}

It turns out that not all moment-angle manifolds are 
of the form \eqref{eq:csum}.  The simplest example of this sort, 
due to Bosio and Meersseman, is presented next. 

\begin{figure}
\setlength{\unitlength}{1cm}
\begin{picture}(4,4)(0,-1.5)
\xygraph{!{0;<1.9cm,0cm>:<0cm,1.9cm>::}
*+{1}="1"-@/^/[rr] *+{2}="2"-@/^/[]!{!(-1,-1.73)}
*+{3}="3"-@/^/"1"-[]!{!(1,-0.15)} *+{4}="4"-[]!{!(0.45,-0.7)}
*+{5}="5"-[]!{!(-0.9,0)}*+{6}="6"-"3"-"5"-"2"-"4"-"6"-"1"
[]!{!(1,-0.58)}*+{7}="7"-"4"
"7"-"6"
"7"-"5"
}
\end{picture}
\caption{\textsf{A triangulation $K$ of  $S^2$ for 
which $\Mo_K$ is not a connected sum of 
products of spheres}}
\label{fig:75}
\end{figure}

\begin{example}
\label{ex:75}
Let $K$ be the triangulation obtained by adding a vertex to 
the face of an octahedron, as shown in Figure~\ref{fig:75}.
The non-edges in $K$ are
$15$, $17$, $26$, $27$, $34$, and $37$.  The 
reader may check that the only pairs of non-edges 
whose union forms a cycle, then, are any choice of 
two of $\set{15, 26, 34}$.   The manifold $\Mo_K$ 
has bigraded Betti numbers as indicated in the 
following diagram, produced by the software 
package Macaulay 2 of Grayson and Stillman \cite{M2}:
\begin{Verbatim}[samepage=true]
          0: 1 .  . . .
          1: . .  . . .
          2: . .  . . .
          3: . 6  . . .
          4: . .  6 . .
          5: . 1  . 1 .
          6: . .  6 . .
          7: . .  . 6 .
          8: . .  . . .
          9: . .  . . .
         10: . .  . . 1
         
\end{Verbatim}

If we denote by $\alpha,\beta$, and $\gamma$, respectively, 
cohomology classes in $H^3(\Mo_K)$ supported on edges 
$15$, $26$, and $34$ in the sense of \eqref{eq:Phimap}, 
the products $\alpha\beta$, $\beta\gamma$, and 
$\gamma\alpha$ are each nonzero, while all 
other products of degree-$3$ classes are zero.  
In \cite[Example 11.5]{BM04}, Bosio and Meersseman 
observe that this manifold $\Mo_K$ cannot be even 
homotopy-equivalent to a manifold of the form 
\eqref{eq:csum}, since the cohomology ring of 
the latter will not have nonzero cup products of 
this form.
\end{example}

\section{Non-formal moment-angle manifolds}
\label{sec:nonformalman}

We now turn to moment-angle manifolds that carry 
non-vanishing Massey products.  We will work 
over a field $\k$ of characteristic $0$.

\subsection{An $8$-vertex triangulation of $S^2$}
\label{subsec:nonformal}

Since the graphs from Figure~\ref{fig:exclude} are planar, 
they can be completed to triangulations of the $2$-sphere.  
These triangulations can be used to produce non-formal 
moment-angle manifolds. 

\begin{example}
\label{ex:eight}
Let $K$ be the flag complex depicted in Figure~\ref{fig:87}.  
The corresponding $11$-dimensional manifold $\Mo_K$ 
has bigraded Betti numbers as listed below.
\begin{Verbatim}[samepage=true]
          0: 1  .  .  .  . .
          1: .  .  .  .  . .
          2: .  .  .  .  . .
          3: . 10  .  .  . .
          4: .  . 16  .  . .
          5: .  .  .  5  . .
          6: .  .  5  .  . .
          7: .  .  . 16  . .
          8: .  .  .  . 10 .
          9: .  .  .  .  . .
         10: .  .  .  .  . .
         11: .  .  .  .  . 1

\end{Verbatim}

By Theorem \ref{th:char}, there are classes 
$\alpha, \beta, \gamma\in H^3(\Mo_K,\k)$ such that 
$\angl{\alpha, \beta, \gamma}\in H^{8}(\Mo_K,\k)$ is 
a non-trivial (decomposable) Massey product, 
with zero indeterminacy. 
\end{example}

\begin{figure}
\setlength{\unitlength}{1cm}
\begin{picture}(5,3.6)(-0.5,-3.3)
\xygraph{
*+{4}="4"-[d]*+{6}="6"-[d]*+{3}="3"
-[r]*+{5}="5"-[u]*+{2}="2"-"6"
-[]!{!(-0.8,-1.4)}*+{1}="1"-@/_/"5"-[]!{!(1.5,0)}*+{7}="7"-"2"
"7"-@/_/"4"
"7"-@/^/"1"-@/^/"4"-"2"
"1"-"3"
"3"-[]!{!(0.5,0.5)}*+{8}="8"-"6"
"2"-"8"-"3"
"8"-"5"
}
\end{picture}
\caption{\textsf{A triangulation $K$ of $S^2$ 
for which $\Mo_K$ is not formal}}
\label{fig:87}
\end{figure}

\begin{remark}
\label{rem:bask}
The manifold $\Mo_K$ from Example \ref{ex:eight} 
belongs to an infinite family of non-formal moment-angle 
manifolds constructed by Baskakov \cite{Ba03}.  This 
is the only member of that family for which $K$ is a 
triangulation of $S^2$. 
\end{remark} 

\begin{remark}
\label{rem:minimal}
This example is minimal: as noted in 
Proposition \ref{prop:formalman}, all triangulations 
of $S^2$ on less than $8$ vertices 
yield formal moment-angle manifolds.  Moreover, 
of the $14$ distinct $8$-point triangulations of $S^2$ 
(as generated by McKay's software \textsf{plantri} \cite{Mac}), 
the flag complex from Example \ref{ex:eight} is the 
unique triangulation $K$ for which $\Mo_K$ possesses 
non-trivial Massey products. 
\end{remark}

\subsection{An infinite family}
\label{subsec:s2}

The next result will permit us to construct many more 
examples of non-formal moment-angle manifolds, not 
covered by Baskakov's method.

\begin{proposition}
\label{prop:cornercut}
If $\Mo_K$ has a non-trivial triple Massey product, 
then so does $\Mo_{K*_F w}$, for any maximal face 
$F$ of dimension at least $2$.

Conversely, if $K$ is a triangulation of $S^2$, $F$ is a triangle, 
and $\Mo_{K*_F w}$ has a non-trivial Massey triple product,
then so does $\Mo_K$.
\end{proposition}

\begin{proof}
By assuming $F$ has dimension at least two, the one-skeleton 
of $K_{123456}$ is unchanged.  By the work above, the non-triviality 
of the Massey product depends only on this one-skeleton.

To prove the converse, it is enough to show that $w$ 
cannot be a vertex $\set{1,\ldots,6}$ of an excluded 
subcomplex as in Figure \ref{fig:excl};  then $K$ and 
$K*_F w$ will have the same subcomplex on 
$\set{1,\ldots,6}$, hence the same non-trivial 
Massey product.  

Consider completing any of the five graphs in 
Figure~\ref{fig:exclude} to a triangulation.  Any vertex on a 
face with more than three edges, including the face at infinity, 
must be joined to a vertex that subdivides that face.  So, 
by inspection, each of $1,\ldots,6$ must have
degree at least $4$ in the $1$-skeleton of any triangulation 
of $S^2$.  By hypothesis, $w$ has degree $3$, so $w\neq i$ 
for each $i$.  This completes the proof.
\end{proof}

Using Example \ref{ex:eight} and Proposition \ref{prop:cornercut}, 
we obtain the following.  

\begin{corollary}
\label{cor:manymassey}
There are infinitely many triangulations $K$ of $S^2$  
for which the moment-angle manifold $\Mo_K$ has non-trivial 
triple Massey products, and thus is not formal. 
\end{corollary} 

\subsection{Asymptotics}
\label{subsec:asym}

It seems likely that ``most'' moment-angle complexes 
are not formal, given a suitable way of making this statement 
precise.  The following result is an example, strengthening
Corollary~\ref{cor:manymassey}.  

\begin{theorem}
\label{th:asym} 
Let $a_n$ denote the number of $n$-point (labeled) 
triangulations of $S^2$, and let $b_n$ denote the number 
of such triangulations $K$ for which the moment-angle 
manifold $\Mo_K$ is formal.  Then 
\[
\lim_{n\to \infty} b_n/a_n=0.
\]
\end{theorem}

\begin{proof}
Let $\mathcal C$ be the set of labeled triangulations of $S^2$.  
Then $\mathcal C$ is small and addable in the terminology of 
McDiarmid, Steger, and Welsh~\cite{MSW05}.  Their Theorem~4.1 
then implies that the probability that a triangulation with $n$ vertices 
chosen uniformly at random contains a fixed, induced subgraph $H$ 
is bounded below by $1-e^{-cn}$, for a positive constant $c$, and for 
sufficiently large $n$.  

Now choose $H$ to be one of the graphs in Figure~\ref{fig:exclude}.  
Using Theorem~\ref{th:char}, we conclude that $\Mo_K$ is not formal, 
for almost all $K\in\mathcal C$.
\end{proof}

We do not know whether the analogous statement 
holds for isomorphism classes of triangulations of $S^2$.

\subsection{Complex moment-angle manifolds}
\label{subsec:complex}

As mentioned in \S\ref{intro:lmv}, work of 
L\'opez de Medrano--Verjovsky \cite{LV97} and 
Meersseman \cite{Me00} makes it possible to construct 
from combinatorial data compact, complex manifolds 
which are not algebraic, and even not K\"ahler or 
symplectic, except in very special cases.

If $K$ is an $n$-vertex polytopal triangulation of $S^m$---that 
is, there exist $n$ points in $\R^{m+1}$ whose convex 
hull is isomorphic to $K$ as a simplicial complex---then 
from \cite[Theorem 12.2]{BM04}, the moment-angle 
manifold $\Mo_K$ admits a complex structure if its 
dimension $n+m+1$ is even.  If $n+m+1$ is odd, then, by 
Remark \ref{rem:extend}, $\Mo_{K,n+1}=\Mo_K\times S^1$, 
and the theorem states that $\Mo_{K,n+1}$ admits a 
complex structure.

For $m\ge 3$, not all triangulations of $S^m$ are polytopal.  
However, a classical theorem of Steinitz states that 
triangulations $K$ of $S^2$ are always polytopal.  
Then our constructions above of triangulations
of $S^2$ for which $\Mo_K$ is not formal apply 
to the corresponding LVM manifolds, too. 

\begin{corollary}
\label{cor:lvm}
There are infinitely many compact, complex manifolds 
arising from the L\'opez de Medrano--Verjovsky--Meersseman 
construction which are not formal. 
\end{corollary}

In particular, by \cite{DGMS}, none of these LVM manifolds 
is K\"{a}hler. 

\subsection{Indecomposable Massey products}
\label{subsec=indecomp}

In view of the theory developed so far, the question 
arises:  are all Massey products in moment-angle 
complexes decomposable?   The next example 
dispels this notion.

\begin{example}
\label{ex:indec}
Consider the $8$-vertex simplicial complex $K$
with Stanley-Reisner ideal 
\begin{equation}
\label{eq:ik}
I_K=(x_1x_2,\, x_1x_3x_4x_5,\, x_3x_4x_5x_6,\, 
x_3x_5x_6x_7,\, x_7x_8).
\end{equation}
(The construction of $K$ is inspired by an example 
of Backelin, reported in \cite{Av98}.)   This simplicial 
complex has $1$ three-face and $12$ four-faces.  
The corresponding $12$-dimensional moment-angle 
complex, $\Mo_K$, has bigraded Betti numbers as 
listed below.

\begin{Verbatim}[samepage=true]
          0: 1 . . . .
          1: . . . . .
          2: . . . . .
          3: . 2 . . .
          4: . . . . .
          5: . . . . .
          6: . . 1 . .
          7: . 3 . . .
          8: . . 4 . .
          9: . . . . .
         10: . . 2 . .
         11: . . . 4 .
         12: . . . . 1

\end{Verbatim}

\begin{claim}
\label{claim:indec}
There is an indecomposable Massey product 
$\angl{\alpha, \beta, \gamma}$ in $H^{12}(\Mo_K,\k)$, 
with $\alpha,\gamma\in H^3(\Mo_K,\k)$ and $\beta\in H^7(\Mo_K,\k)$.
\end{claim}

\begin{proof}
Label the generators of $I_K$ as $e_1,\dots ,e_5$, and 
make the Taylor resolution on the exterior algebra generated by 
$e_1,\dots ,e_5$ (see e.g.~\cite{Yu99}).  A computation 
shows that  $e_1 e_3 = - d(e_1e_2e_3)$, while 
$e_3 e_5 = - d(e_3e_4e_5)$.  Thus, the Massey product 
$\angl{e_1,e_3,e_5}$ is defined.  

Note that, for each $p\ne 6$, either $b_{p}(\Mo_K)=0$ or 
$b_{12-p}(\Mo_K)=0$, and so $H^{p}\cdot H^{12-p}=0$.  
Furthermore, $b_6(\Mo_K)=1$, 
and so, by formula \eqref{eq:baskakov}, we have 
$H^6\cdot H^6=0$ as well. Thus, 
$\angl{e_1,e_3,e_5}$ has zero indeterminacy, 
and is indecomposable. 
\end{proof}

Following Backelin~\cite{Bac82} as in Remark~\ref{rem:golod}, 
it is an exercise to show here that
\begin{equation}
\label{eq:back}
\sum_{p,q\geq0}\dim_\k\Tor^{S/I}_{p,q}(\k,\k)s^qt^p
= \frac{(1+st)^6}{(1-st-3s^6t^2-s^7t^3)(1-st)},
\end{equation}
which from \eqref{eq:hstensor}, means
\begin{equation}
\label{eq:hug}
h(U(\g_\Mo),t)=((1-t-3t^6-t^7)(1-t)(1+t)^2)^{-1}.
\end{equation}

The identity \eqref{eq:ranks} makes it possible to compute the ranks
$\phi_r$ of $\pi_r(\Mo_K)$ recursively, and we find 
\begin{equation}
\label{eq:piranks}
\begin{array}{c|rrrrrrrrrrrrr}
r & 1 & 2 & 3 & 4 & 5 & 6 & 7 & 8 & 9 & 10 & 11 & 12 & 13\\ \hline
\phi_r & 0 & 0& 2 &0 &0 &0 &3 &4 &4 &4 &4 & 4 & 7\\ 
\end{array}
\end{equation}

The indecomposable Massey product above contributes 
a nonzero differential in the $E_2$ term of the 
Eilenberg-Moore spectral sequence, so the
Lie algebra of primitives in  $\Ext_{H^*(\Mo_K,\k)}(\k,\k)$ 
agrees with $\pi_*(\Omega\Mo_K)\otimes \k$ up to but not including 
$\pi_{10}(\Omega\Mo_K)\otimes \k$.
\end{example}

\subsection{Bier spheres}
\label{subsec:bier}

We conclude this section with an example of a moment-angle 
manifold whose cohomology has indecomposable triple 
Massey products. 
Our approach uses a construction of triangulated 
spheres, due to Bier \cite{Bier}.  We briefly describe 
this construction here, and refer to \cite{Ma03}, \cite{dL04} for
a complete discussion. 

Let $K$ be proper simplicial subcomplex of the $n$-simplex.  
Then $K^{\star}$, its {\em combinatorial Alexander dual}, is 
defined to be the simplicial complex whose simplices are the
complements in $[n]$ of the non-simplices of $K$.  
The {\em Bier sphere} associated to $K$ is the deleted
join of $K$ with $K^{\star}$:
\begin{equation}
\label{eq:bier}
\Bier(K) = \{ \sigma * \tau \in K * K^{\star} 
\mid \sigma \cap \tau =\emptyset \}.
\end{equation}

Remarkably, $\widehat{K}=\Bier(K)$ is a triangulation of 
$S^{n-2}$ on $2n$ vertices: see  \cite{Ma03}, \cite{dL04}.  
Consequently, the moment-angle complex $\Mo_{\widehat{K}}$ 
is a manifold, of dimension $3n-1$. 

\begin{example}
\label{ex:simplebier}
If $K$ is the $4$-gon, then $K^{\star}$ is the disjoint union 
of two edges, $\widehat{K}=\Bier(K)$ is the $8$-vertex 
triangulation of $S^2$ from Figure \ref{fig:87}, and 
$\Mo^{11}_{\widehat{K}}$ is the non-formal manifold 
from Example \ref{ex:eight}. 
\end{example}

\begin{example}
\label{ex:indecbier}
Let $K$ be the $8$-vertex simplicial complex from 
Example \ref{ex:indec}, and let $\widehat{K}=\Bier(K)$, 
a $16$-vertex triangulation of $S^{6}$.  The 
manifold $\Mo^{23}_{\widehat{K}}$ has 
Betti numbers 
$1,0,0,10,4,12,98,130,91,233,377,268,
268,377,233,91,130,98,12,4,10,0,0,1$.

\begin{claim}
\label{claim:indecbier}
The manifold $\Mo_{\widehat{K}}$ carries an 
indecomposable Massey triple product 
$\angl{\alpha, \beta, \gamma}\in H^{12}(\Mo_{\widehat{K}}, \k)$.  
\end{claim}

\begin{proof}
Notice that $K$ is a full subcomplex of $\widehat{K}$, 
so we may apply the results of \S\ref{subsec:cellular}.  Using
Lemma~\ref{lem:cellmas}, the Massey product in 
$C^*_{\CW}(\Mo_K,\k)$ lifts through the split surjection
$\Mo_f^* \colon C^*_{\CW}(\Mo_{\widehat{K}},\k)  \to 
C^*_{\CW}(\Mo_K,\k)$ to a Massey product for 
$\Mo_{\widehat{K}}$.
\end{proof}

\end{example}

\section{Subspace arrangements}
\label{sec:subspaces}

Let $\A$ be an arrangement of linear subspaces in 
$\C^n$.  The {\em intersection lattice} of the arrangement, 
$L(\A)$, is the poset of all intersections among the 
subspaces in $\A$, ordered by reverse inclusion.  
The {\em complement} of the arrangement is simply the 
space $X(\A)=\C^n\setminus\bigcup_{H\in \A} H$.  
A fundamental question in the subject is whether the 
homotopy type of the complement is completely determined 
by the combinatorial data, i.e., the intersection lattice and the 
codimensions of the subspaces.

\subsection{Cohomology}
\label{subsec:sscoho}

In \cite{GM88}, Goresky and MacPherson gave a combinatorial 
formula for the Betti numbers of $X(\A)$, one version of which
we recall here:
\begin{equation}
\label{eq:gm}
H^i(X(A),\k) = \bigoplus_{S\in L(\A)} \widetilde{H}_{2\rho(S)-i-2}(
(\widehat{0},S),\k),
\end{equation}
where $(\widehat{0},S)$ denotes an open interval in the 
order complex of $L(\A)$, and $\rho(S)$ is the complex 
codimension of $S$ in $\C^n$.  
Later on, the cohomology ring $H^*(X(\A),\Z)$ was 
computed (also in purely combinatorial terms) by 
de Longueville and Schultz \cite{dLS01} and 
Deligne, Goresky, and MacPherson \cite{DGM00}. 

De Concini and Procesi observe in \cite{dCP95} that the 
mixed Hodge structure on the cohomology of a complex 
subspace complement of weight $2r$ has type $(r,r)$, 
for each $r$.  More precisely, Deligne, Goresky
and MacPherson note \cite[Exemple 1.14]{DGM00}
that the summand of \eqref{eq:gm} labelled by $S$ is 
of type $(r,r)$, where $r=2\rho(S)$.

\subsection{Coordinate subspaces}
\label{subsec:sscoord}

Now let $K$ be a finite simplicial complex on vertex set $[n]$.
Associated to $K$ there is a subspace arrangement 
$\A_{K}$ in $\C^n$, consisting of a coordinate subspace 
$H_{\sigma}=\{ z\in \C^n \mid z_i =0 \ \text{if $i\notin \sigma$}\}$ 
for each non-empty simplex $\sigma\subseteq K$.  
It is readily seen that $X(\A_{K})=\Mo_K(\C,\C^*)$; thus, 
$X(\A_{K})$ deform-retracts onto the moment-angle 
complex $\Mo_K=\Mo_K(D^2,S^1)$, cf.~\cite{BP00}.

The subspaces in $\A_K$ are defined by equations
$\set{z_{i_1}=z_{i_2}=\cdots=z_{i_k}=0}$, over index 
sets for which $x_{i_1}\cdots x_{i_k}$ is a generator 
of the Stanley-Reisner ideal $I=I_K$.
Along the same lines, the intersection lattice $L(\A)$ 
is isomorphic to both the lattice of least common 
multiples of the generators of $I$ (studied in \cite{GPW99}), 
and the opposite to the face lattice of the Alexander dual 
simplicial complex $K^\star$.

Using the second interpretation, Buchstaber and Panov 
observe in \cite{BP00} that the Goresky-MacPherson 
formula \eqref{eq:gm} specializes for coordinate 
subspace arrangements to the classical Hochster 
formula \eqref{eq:hochster}, via Alexander 
duality.  Indeed, for $p,q\geq0$,
\begin{eqnarray}
\label{eq:hochsteragain}
\Tor^S_{p-q,p+q}(S/I,\k) &\cong & \bigoplus_{\abs{\I}=p}
\widetilde{H}^{q-1}(K_\I,\k)\\
&\cong & \bigoplus_{\abs{\I}=p} 
\widetilde{H}_{p-q-2}((K_\I)^\star,\k)\notag \\
&\cong& \bigoplus_{\stackrel{\sigma\in K^\star}{\abs{\sigma}=n-p}}
\widetilde{H}_{p-q-2}({\rm link}_{K^\star}(\sigma),\k).\notag
\end{eqnarray}
Since $H^i(X(\A_K),\k)=H^i(\Mo_K,\k)=\Tor^S_{*,i}(S/I,\k)$,
the argument is completed by noting that the link of a simplex 
$\sigma$ in $K^\star$ is homeomorphic to its barycentric subdivision, 
the interval $(\widehat{0},[n]-\sigma)$ in $L(\A)$.  We refer to
Eagon and Reiner \cite{ER98} for a related discussion. 

The cohomology of the subspace complement $X(\A_K)$ inherits 
a bigrading from $\Tor$.   Let $\gr^W$ denote the grading 
associated to the weight filtration.  Under the identification 
$H^q(X(\A_K),\C)=\Tor^S_{*,q}(S/I,\C)$, then,
keeping track of the bigrading above gives
\begin{equation}
\label{eq:weight}
\gr^W_r H^*(X(\A_K),\C) = \bigoplus_{p}\Tor^S_{p,2r-p}(S/I,\C),
\end{equation}
for all $r\geq0$.

\subsection{The formality question}
\label{ss:formality}

In \cite{dCP95} de~Concini and Procesi gave a DGA model 
for the rational cohomology ring of a complex subspace 
arrangement $\A$. Using this model (as simplified by 
Yuzvinsky in \cite{Yu02}), Feichtner and Yuzvinsky \cite{FY05} 
prove the following:  If $L(\A)$ is a geometric lattice, then 
$X(\A)$ is a formal space.  

Among arrangements with geometric intersection  
lattice, best understood are the ``redundant" 
subspace arrangements, for which explicit 
computations of the homotopy Lie algebra can 
be done, under some assumptions; 
see \cite{CCX03, PS04, DS06}.
However, it follows from \cite{Ba03} and the work above that the 
complement of a subspace arrangement need {\em not}\/ 
be formal in general.  Indeed, we have the following result. 

\begin{proposition}
\label{prop:non formal subspace arr}
Suppose $K$ is a simplicial complex, containing as a full 
subcomplex one of the graphs listed in Figure~\ref{fig:exclude}. 
Then the complement of the coordinate subspace arrangement 
$\A_K$ is not formal.
\end{proposition}

\begin{proof}
Recall the complement $X(\A_K)$ is homotopy equivalent 
to the moment-angle complex $\Mo_K$.  The result thus follows
from Theorem~\ref{th:char}.
\end{proof}

The simplest such example is obtained by letting $K$ be the simplicial
complex of Figure \ref{fig:sixvertex}.  This is a simplicial complex
on $6$ vertices, with five minimal non-faces, and so $\A_K=\{H_1,\dots
,H_5\}$ where $H_i= \set{z \in \C^6 \mid z_i=z_{i+1}=0}$.  As
expected, the intersection lattice of $\A_K$ is not geometric, either:
the fragment of the lattice depicted in Figure \ref{fig:lattice} makes
it clear why.

\begin{figure}
\setlength{\unitlength}{1cm}
\begin{picture}(4,2.9)(0.9,0)
\xygraph{
*+{12}="12"-[ur]*+{123}="123"
(
,-[dr]*+{13}="13"-[ur]*+{134}="134"-[dr]*+{34}="34"
,-[ur]*+{1234}-[dr]*+{134}
)
}
\end{picture}
\caption{\textsf{A non-geometric lattice}}
\label{fig:lattice}
\end{figure}

\begin{ack}
The first author would like to thank Srikanth Iyengar for helpful
conversations regarding Golod maps and formality in commutative
algebra.  Both authors thank the referee for helpful suggestions 
that led to improvements in the presentation and content of the paper.
\end{ack}

\providecommand{\bysame}{\leavevmode\hbox to3em{\hrulefill}\thinspace}
\providecommand{\MR}{\relax\ifhmode\unskip\space\fi MR }
\providecommand{\MRhref}[2]{%
  \href{http://www.ams.org/mathscinet-getitem?mr=#1}{#2}
}
\providecommand{\href}[2]{#2}
\newcommand{\arxiv}[1]{{\texttt{\href{http://arxiv.org/abs/#1}%
{arXiv:{#1}}}}}

\end{document}